\documentclass[11pt]{amsart}
\usepackage{amssymb}
\usepackage{amsthm}

\usepackage{amsfonts}
\usepackage{amsmath}
\usepackage[abbrev]{amsrefs}
\usepackage{enumerate}

\usepackage{hyperref}

\newtheorem{Theorem}{Theorem}[section]

\newtheorem{Lemma}[Theorem]{Lemma}
\newtheorem{Corollary}[Theorem]{Corollary}
\newtheorem{Proposition}[Theorem]{Proposition}

\theoremstyle{definition}

\newcommand{\R}{\mathbb{R}}

\newcommand{\N}{\mathbb{N}}

\DeclareMathOperator{\esssup}{ess\,sup}

\DeclareMathOperator*{\sign}{sign}
\DeclareMathOperator{\supp}{supp}

\begin{document}

\title{The Daugavet property in the Musielak-Orlicz spaces}

\author{Anna Kami\'{n}ska}
\address{Department of Mathematical Sciences,
The University of Memphis, TN 38152-3240}
\email{kaminska@memphis.edu}

\author{Damian Kubiak}
\address{Mathematics Department,
Tennessee Technological University, 110 University Drive,
Box 5054, Cookeville, TN 38505}
\email{dkubiak@tntech.edu}
\thanks{The second author was supported by the Tennessee Technological University internal Faculty Research Grant during 2013-14.}
\thanks{\underline{Please cite this article in press as:} A. Kami\'{n}ska, D. Kubiak, The Daugavet property in the Musielak-Orlicz spaces, J. Math. Anal. Appl. 427 (2) (2015), 873-898, \href{http://www.ams.org/mathscinet-getitem?mr=3323013}{MR3323013}, \url{http://dx.doi.org/10.1016/j.jmaa.2015.02.035}.}

\date{October 13, 2014; Revised March 03, 2015}

\begin{abstract}
We show that among all Musielak-Orlicz function spaces on a $\sigma$-finite non-atomic complete measure space equipped with either the Luxemburg norm or the Orlicz norm the only spaces with the Daugavet property are $L_1$, $L_{\infty}$, $L_1\oplus_1 L_{\infty}$ and $L_1\oplus_{\infty} L_{\infty}$. In particular, we obtain  complete characterizations of the Daugavet property in the weighted interpolation spaces, the variable exponent Lebesgue spaces (Nakano spaces) and the Orlicz spaces. 
\end{abstract}

\maketitle

\section{Introduction}

Let $(X,\|\cdot\|)$ be a Banach space and $T:X\to X$ be a bounded linear operator. The equation 
\begin{equation}
\label{DE}
\|T+I\|=1+\|T\|\text{,}
\end{equation}
where $I$ is the identity operator on $X$, is called the Daugavet equation. If the Daugavet equation is satisfied by every rank one operator $T$ then $X$ is said to have the Daugavet property. It is known that if $X$ has the Daugavet property then Eq. (\ref{DE}) is satisfied by every weakly compact operator.

Eq. (\ref{DE}) was first studied by I.K. Daugavet in the space $C(0,1)$ \cite{MR0157225}. Examples of spaces which have the Daugavet property are $L_1$ and $L_{\infty}$ over a non-atomic measure space as well as $C(K)$, where $K$ is a compact Hausdorff space with no isolated points. 
Moreover, finite direct sums $\oplus_1$ and $\oplus_{\infty}$ of spaces with the Daugavet property possess that property as well \cite{MR1126202}. It is known that Banach spaces with the Daugavet property fail the Radon-Nikodym property and do not embed into a space with an unconditional basis. For a historical overview on the Daugavet property we refer to \cite{MR3013764}. An introductory exposition on the Daugavet equation and the Daugavet property can be found in \cite{MR1921782}.

It has been recently showed that among all rearrangement invariant function spaces over a non-atomic finite measure spaces only $L_1$ and $L_{\infty}$ have the Daugavet property \cite{AKM_Dri, MR3028566}. Inspired by that result we study the Daugavet property in the class of Musielak-Orlicz function spaces on a $\sigma$-finite non-atomic complete measure space. These spaces are not rearrangement invariant in general. The variable exponent Lebesgue spaces (Nakano spaces) and the Orlicz spaces appear as special cases of the Musielak-Orlicz spaces. It should be mentioned that the class of Musielak-Orlicz spaces we consider here is the most general and includes also the interpolation spaces $L_1\cap L_{\infty}$ and $L_1+L_{\infty}$, as well as their weighted versions which are studied in Section 3. In section 4, using an observation that the unit sphere of a Banach space with the Daugavet property does not contain a uniformly non-$\ell^2_1$ point, we prove that the only Musielak-Orlicz spaces (equipped with either the Luxemburg norm or the Orlicz norm) with the Daugavet property are $L_1$, $L_{\infty}$, $L_1\oplus_1 L_{\infty}$ or $L_1\oplus_{\infty} L_{\infty}$ with weights. This generalizes several results obtained earlier in \cite{MR3013764} and solves the problem of the Daugavet property in the Musielak-Orlicz spaces completely. In the appendix we give a proof of K\"{o}the duality in the most general case of the Musielak-Orlicz spaces. That result is of course well known but it seems that a direct proof in such generality has never been published. 

\section{Preliminaries}

For a Banach space $X$, by $S(X)$ and $B(X)$ we denote the unit sphere and the unit ball, respectively. The space of all bounded linear functionals on $X$ is denoted by $X^*$. 
Let $(X,\|\cdot\|)$ be a (real) Banach space. 
For any $x^{*}\in S(X^{*})$ and  $\epsilon>0$ the set 
\begin{equation*}
S(x^{*};\epsilon)=\{x\in B(X):x^{*}x>1-\epsilon\}
\end{equation*} is called \emph{a slice} determined by $x^{*}$ and $\epsilon$. We say that a Banach space $(X,\|\cdot\|)$ has the \emph{slice (or local) diameter $2$ property} if every slice of $B(X)$ has diameter $2$. It is known that every space with the Daugavet property has the slice diameter $2$ property \cite{MR3098474}.

Let $X$ and $Y$ be Banach spaces with norms $\|\cdot\|_X$ and $\|\cdot\|_Y$, respectively. By $X\oplus_1 Y$ we denote the Banach space consisting of all ordered pairs $(x,y)$ where $x\in X$ and $y\in Y$ with the norm $\|(x,y)\|=\|x\|_X+\|y\|_Y$. Similarly, by $X\oplus_{\infty} Y$ we denote the Banach space consisting of all ordered pairs $(x,y)$ where $x\in X$ and $y\in Y$ with the norm $\|(x,y)\|=\max\{\|x\|_X,\|y\|_Y\}$. It is clear that $(X\oplus_1 Y)^*=X^*\oplus_{\infty} Y^*$ and $(X\oplus_{\infty} Y)^*=X^*\oplus_1 Y^*$ with equality of norms. The fact that two Banach spaces $X$ and $Y$ are isometrically isomorphic is denoted by $X\simeq Y$ or by $X=Y$ if an isometric isomorphism between $X$ and $Y$ is obvious (for example the identity mapping).

In the sequel we assume that $(\Omega,\Sigma,\mu)$ is a non-atomic $\sigma$-finite complete measure space. By $L_0=L_0(\Omega)$ we denote the set of all (equivalence classes with respect to the equality $\mu$-a.e. of) measurable extended-real valued functions on $\Omega$.

Let $(X,\|\cdot\|_X)$ be a Banach function lattice on $(\Omega,\Sigma,\mu)$, that is $X\subset L_0$ and if $|x|\leqslant |y|$ $\mu$-a.e. on $\Omega$, $x\in L_0$, $y\in X$ then $x\in X$ and $\|x\|_X\leqslant\|y\|_X$. For any function $x\in L_0$, the support of $x$ is defined by $\supp(x)=\{t\in\Omega:x(t)\ne 0\}$. Recall that $\supp(X)$ is a measurable subset of $\Omega$ such that every element of $X$ vanishes $\mu$-a.e. on $\Omega\setminus\supp(X)$ and for every measurable subset $E$ of $\supp(X)$ with positive measure there is a measurable set $F\subset E$ with finite and positive measure such that $\chi_F\in X$ \cite[p. 14]{MR1996919}. 
An element $x\in X$ is called order continuous if for every $0\leqslant x_n\leqslant |x|$ such that $x_n\downarrow 0$ $\mu$-a.e. it holds $\|x_n\|_X\to 0$. By $X_a$ we denote the set of all order continuous elements of $X$. A Banach function lattice $X$ is said to have the Fatou property whenever for any sequence $(x_n)$ in $X$ and $x\in L_0$ such that $x_n\to x$ $\mu$-a.e. on $\Omega$ and $\sup\|x_n\|_X<\infty$, we have that $x\in X$ and $\|x\|_X\leqslant \liminf \|x_n\|_X$. Given a measurable set $\Gamma\subset\Omega$ we denote $X(\Gamma)=\{x\in X:\mu(\supp(x)\setminus\Gamma)=0\}$ and the norm $\|\cdot\|_{X(\Gamma)}$ on $X(\Gamma)$ is defined by $\|x\|_{X(\Gamma)}=\|x\chi_{\Gamma}\|_X$. It is clear that $X(\Gamma)$ is continuously embedded in $L_0$.

Let $(X,\|\cdot\|_X)$ and $(Y,\|\cdot\|_Y)$ be Banach function lattices on $(\Omega,\Sigma,\mu)$ with the Fatou property and let $\Gamma_1,\Gamma_2$ be measurable sets such that $\Gamma_2\subset\Gamma_1\subset\Omega$. We define 
\begin{equation*}
X(\Gamma_1)\cap Y(\Gamma_2)=Y(\Gamma_2)\cap X(\Gamma_1)=\{x\in L_0:x\in X(\Gamma_1)\text{ and }x\chi_{\Gamma_2}\in Y(\Gamma_2)\}\text{,}
\end{equation*}
and equip it with the norm
\begin{equation*}
\|x\|=\|x\|_{X(\Gamma_1)\cap Y(\Gamma_2)}=\|x\|_{Y(\Gamma_2)\cap X(\Gamma_1)}=\max\{\|x\|_{X(\Gamma_1)},\|x\chi_{\Gamma_2}\|_{Y(\Gamma_2)}\}\text{.}
\end{equation*}
It follows that $(X(\Gamma_1)\cap Y(\Gamma_2),\|\cdot\|)$ is a Banach function lattice on $(\Omega,\Sigma,\mu)$ with the Fatou property and $\supp(X(\Gamma_1)\cap Y(\Gamma_2))=\Gamma_1$. Moreover, we set
\begin{equation*}
X(\Gamma_1)+Y(\Gamma_2)=Y(\Gamma_2)+X(\Gamma_1)=\{z\in L_0:z=x+y\text{ for some }x\in X(\Gamma_1)\text{ and }y\in Y(\Gamma_2)\}\text{,}
\end{equation*}
and equip it with the norm
\begin{equation*}
\|z\|^{\Sigma}=\|z\|^{\Sigma}_{X(\Gamma_1)+Y(\Gamma_2)}=\|z\|^{\Sigma}_{Y(\Gamma_2)+X(\Gamma_1)}=\inf\{\|x\|_{X(\Gamma_1)}+\|y\|_{Y(\Gamma_2)}:z=x+y, x\in X(\Gamma_1), y\in Y(\Gamma_2)\}\text{.}
\end{equation*}
Again, $(X(\Gamma_1)+Y(\Gamma_2),\|\cdot\|^{\Sigma})$ is a Banach function lattice on $(\Omega,\Sigma,\mu)$ with the Fatou property and $\supp(X(\Gamma_1)+Y(\Gamma_2))=\Gamma_1$.
The proof of completeness of $X(\Gamma_1)\cap Y(\Gamma_2)$ and $X(\Gamma_1)+Y(\Gamma_2)$ is essentially the same as the proof of Theorem 1.3 \cite[p. 97]{MR928802}. 

It should be noted that given two Banach function lattices $X$ and $Y$ on $(\Omega,\Sigma,\mu)$ the space $X\cap Y$ is usually defined (explicitly or implicitly) as the set of functions belonging to both, $X$ and $Y$ (see for example \cite[p. 97]{MR928802},\cite[p. 9]{MR649411} or \cite[p. 16]{MR1996919}). However, as we will see in Theorem \ref{MOdecomp}, we need to consider spaces $X\cap Y$ as defined in the previous paragraph.

Recall \cite[p. 44]{MR649411} that the K\"{o}the dual $X^{\prime}$ of $X$ is the collection of those $y\in L_0$ such that $\supp(y)\subset\supp(X)$ and 
\begin{equation*}
\|y\|_{X^{\prime}}=\sup\left\{\int_{\Omega}|xy|d\,\mu :\|x\|_X\leqslant 1\right\}<\infty\text{.}
\end{equation*}
The space $(X^{\prime},\|\cdot\|_{X^{\prime}})$ is a Banach function lattice on $(\Omega,\Sigma,\mu)$ with the Fatou property and $\supp(X^{\prime})=\supp(X)$. A functional $F\in X^*$ is said to be order continuous if $F(x_n)\to 0$ whenever $x_n, x\in X$, $x_n\to 0$ and $|x_n|\leqslant x$ $\mu$-a.e. on $\Omega$. It is known that $F\in X_c^*$, the set of all order continuous functionals on $X$, if and only if there exists a unique $y\in X^{\prime}$ such that $F(x)=\int_{\Omega}xy\,d\mu$ and $\|F\|=\|y\|_{X^{\prime}}$. Hence $X^{\prime}$ is isometrically isomorphic to $X_c^*$. Moreover $X^{\prime\prime}=X$ with equality of norms if and only if $X$ has the Fatou property. For details on Banach lattices see \cite{MR928802, MR649411, MR0511676, MR1996919}.

Recall also the following result on K\"{o}the duality of spaces $X\cap Y$ and $X+Y$ (cf. Lemma 1.12 \cite[p. 18]{MR1996919}).
\begin{Theorem}
\label{dualXY}
Let $(X,\|\cdot\|_X)$ and $(Y,\|\cdot\|_Y)$ be Banach function lattices on $(\Omega,\Sigma,\mu)$ with the Fatou property and $\Gamma_2\subset\Gamma_1\subset\Omega$, where $\Gamma_1$, $\Gamma_2$ are measurable sets of positive measure. The following K\"{o}the dualities hold true,
\begin{equation*}
(X(\Gamma_1)+Y(\Gamma_2))^{\prime}=X(\Gamma_1)^{\prime}\cap Y(\Gamma_2)^{\prime}
\end{equation*}
and
\begin{equation*}
(X(\Gamma_1)\cap Y(\Gamma_2))^{\prime}=X(\Gamma_1)^{\prime}+Y(\Gamma_2)^{\prime}
\end{equation*}
with equality of norms.
\end{Theorem}

A function $\varphi:[0,\infty)\to[0,\infty]$ is called an Orlicz function, if $\varphi$ is not identically $0$, $\lim_{u\to 0^+}\varphi(u)=\varphi(0)=0$, and $\varphi$ is left continuous and convex on $(0,b_{\varphi}]$, where $b_{\varphi}=\sup\{u>0:\varphi(u)<\infty\}$. It follows that $\varphi$ is continuous on $(0,b_\varphi)$. For an Orlicz function $\varphi$ we define $a_{\varphi}=\sup\{u\geqslant 0:\varphi(u)=0\}$ and $d_{\varphi}=\sup\{u\in[0,b_{\varphi}):\varphi(u/2)=\varphi(u)/2\}$. 
Clearly $0\leqslant a_{\varphi}\leqslant d_{\varphi}\leqslant b_{\varphi}\leqslant\infty$, $a_{\varphi}<\infty$, and $b_{\varphi}>0$. Moreover, if $a_{\varphi}=0$ and $d_{\varphi}>0$ then $\varphi(u)=cu$ for all $u\in[0,d_{\varphi})$ and some constant $c>0$. Clearly, if $a_{\varphi}>0$ then $d_{\varphi}=a_{\varphi}$. 
For convenience we denote $\varphi(\infty)=\infty$.

A function $M:\Omega\times[0,\infty)\to[0,\infty]$ is called a Musielak-Orlicz function if for $\mu$-a.e. $t\in \Omega$, $M(t,\cdot)$ is an Orlicz function and for all $u\geqslant 0$, $M(\cdot,u)$ is measurable.

For a Musielak-Orlicz function $M$ we define the functions $a_M(t)=\sup\{u\geqslant 0:M(t,u)=0\}$, $b_M(t)=\sup\{u>0:M(t,u)<\infty\}$ and $d_M(t)=\sup\{u\in[0,b_M(t)):M(t,u/2)=M(t,u)/2\}$. The functions $a_M$, $b_M$ and $d_M$ are measurable \cite{MR1410390, MR2784391}. The basic pointwise properties of the functions $a_M$, $b_M$ and $d_M$ follow from the above discussion on Orlicz functions.

For a Musielak-Orlicz function $M$ the complementary function $N:\Omega\times[0,\infty)\to[0,\infty]$ is defined by
\begin{equation*}
N(t,u)=\sup\{uv-M(t,v):v>0\}\text{, }t\in \Omega\text{ and }u\geqslant 0\text{.}
\end{equation*}
It is known, that the function $N$ complementary to $M$ is a Musielak-Orlicz function.

Let $M$ and $N$ be complementary Musielak-Orlicz functions. The semimodular $\rho_M:L_0\to[0,\infty]$ given by 
\begin{equation*}
\rho_M(x)=\int_{\Omega}M(t,|x(t)|)\, d\mu
\end{equation*} is convex and defines the Musielak-Orlicz function space 
\begin{equation*}
L_M=\{x\in L_0:\rho_M(\lambda x)<\infty\text{ for some }\lambda>0\}
\end{equation*} with the Luxemburg norm 
\begin{equation*}
\|x\|_M=\inf\{\lambda>0:\rho_M(x/\lambda)\leqslant 1\}\text{.}
\end{equation*} 
We also consider the space $E_M$ of all \emph{finite elements} in $L_M$,
\begin{equation*}
E_M=\{x\in L_0:\rho_M(\lambda x)<\infty\text{ for all }\lambda>0\}\text{.}
\end{equation*}
For details on modular spaces and (semi)modulars we refer to \cite{MR724434}.

 Recall that the K\"othe dual $(L_M,\|\cdot\|_M)^{\prime}=(L_N,\|\cdot\|_N^o)$, where 
\begin{equation*}
\|x\|_N^o=\sup\left\{\int_{\Omega} xy\, d\mu:\rho_M(y)\leqslant 1\right\}
\end{equation*} is the Orlicz norm on $L_N$ (see Theorem \ref{dualLM}). The Orlicz norm is equal to the Amemiya norm \cite{MR2286921,MR1909821}, that is 
\begin{equation*}
\|x\|_N^o=\inf_{k>0}\frac{1}{k}[1+\rho_N(kx)]\text{.}
\end{equation*}
Moreover $(L_M,\|\cdot\|_M^o)^{\prime}=(L_N,\|\cdot\|_N)$ (see Theorem \ref{OrliczNormDual}). It is known that $\|x\|_M\leqslant\|x\|_M^o\leqslant2\|x\|_M$ for all $x\in L_M$ \cite[p. 9]{MR724434}.
 In the sequel, we denote $L_M=(L_M,\|\cdot\|_M)$ and $L_N^o=(L_N,\|\cdot\|_N^o)$. It follows from the general theory of Banach lattices that $L_M^*$ is isomorphic to $(L_M)^{\prime}\oplus S$, where $S$ is the set of all singular functionals on $L_M$. Moreover, every singular functional evaluates to $0$ at order continuous elements of $L_M$ \cite{luxemburg1983riesz}. 
Since  $(L_M)^{\prime}=L_N^o$ so $L_N^o$ is isometrically isomorphic to a subspace of $L_M^*$.

If a Musielak-Orlicz function $M(t,u)=\varphi(u)$ for all $t\in\Omega$, where $\varphi$ is an Orlicz function, then $L_M=L_{\varphi}$ is called an Orlicz space. In this case we denote $\psi(u)=N(t,u)$, where $N$ is the complementary function of $M$ \cite{MR1410390}.

Let $p\in L_0$ be such that $1\leqslant p(t)\leqslant\infty$. If, for $\mu$-a.e. $t\in \Omega$, 
\begin{equation*}
M(t,u)=\begin{cases}
\frac{u^{p(t)}}{p(t)}\text{ if }t\in \Omega\setminus\Omega_{\infty}\text{,}\\
\alpha(u)\text{ if }t\in \Omega_{\infty}\text{,}
\end{cases}
\end{equation*}
where $\Omega_{\infty}=\{t\in \Omega:p(t)=\infty\}$ and $\alpha(u)=0$ for $0\leqslant u\leqslant 1$ and $\alpha(u)=\infty$ for $u>1$ then the space $L_M$ is called a variable exponent Lebesgue space (or Nakano space) and is denoted by $L_{p(t)}$.

We use the standard convention that $a/0=\infty$, $a/\infty=0$, $\infty/a=\infty$ for $a\in(0,\infty)$, $0\cdot\infty=0$, $\infty\leqslant\infty$, $\inf\emptyset=\infty$.

\section{Weighted Interpolation Spaces $L_{1,v}+L_{\infty,w}$ and $L_{1,w}\cap L_{\infty,v}$}

In this section we study the Daugavet property in weighted interpolation spaces $L_{1,v}+L_{\infty,w}$ and $L_{1,w}\cap L_{\infty,v}$, which are in fact, as we will see in the next section, the Musielak-Orlicz spaces generated by certain Musielak-Orlicz functions. We prove the criteria of the Daugavet property in both spaces which will be applied in the proofs of main results in section 4. 

For a measurable set $\Gamma\subset\Omega$, a function $u\in L_0(\Omega)$ is called a weight function on $\Gamma$ if $0<u<\infty$ $\mu$-a.e. on $\Gamma$. Given an arbitrary weight function $u$ on $\Omega$, we denote 
\begin{equation*}
L_{\infty,u}=L_{\infty,u}(\Omega)=\{x\in L_0:xu\in L_{\infty}(\Omega)\}
\end{equation*}
and equip it with the standard norm 
\begin{equation*}
\|x\|_{\infty,u}=\|xu\|_{\infty}=\esssup_{t\in \Omega}|x(t)u(t)|\text{.}
\end{equation*}
Similarly
\begin{equation*}
L_{1,u}=L_{1,u}(\Omega)=\{x\in L_0:xu\in L_1(\Omega)\}
\end{equation*}
and 
\begin{equation*}
\|x\|_{1,u}=\|xu\|_{1}=\int_{\Omega}|xu|\, d\mu\text{.}
\end{equation*}

Let $\Gamma\subset\Omega$ be a measurable set with $\mu(\Gamma)>0$ and $u$ be a weight function on $\Gamma$. It follows that 
\begin{equation*}
L_{\infty,u}(\Gamma)=\{x\in L_0:\mu(\supp(x)\setminus\Gamma)=0\text{ and }xu\in L_{\infty}(\Gamma)\}
\end{equation*}
and for $x\in L_{\infty,u}(\Gamma)$, $\|x\|_{\infty,u}=\esssup_{t\in \Gamma}|x(t)u(t)|$.
Similarly 
\begin{equation*}
L_{1,u}(\Gamma)=\{x\in L_0:\mu(\supp(x)\setminus\Gamma)=0\text{ and }xu\in L_1(\Gamma)\}
\end{equation*}
and for $x\in L_{1,u}(\Gamma)$, $\|x\|_{1,u}=\int_{\Gamma}|xu|\, d\mu$.
Both spaces $(L_{\infty,u}(\Gamma),\|\cdot\|_{\infty,u})$ and $(L_{1,u}(\Gamma),\|\cdot\|_{1,u})$ are Banach function lattices on $(\Omega,\Sigma,\mu)$ with the Fatou property continuously embedded in $L_0$ and with $\supp(L_{\infty,u}(\Gamma))=\supp(L_{1,u}(\Gamma))=\Gamma$. 

 It is clear that $L_{1,u}(\Gamma)\simeq L_1(\Gamma)$ by the isometric isomorphism from $L_{1,u}(\Gamma)$ to $L_1(\Gamma)$ mapping $x\mapsto xu$. Similarly, $L_{\infty,u}(\Gamma)\simeq L_{\infty}(\Gamma)$. It follows that $(L_{1,u}(\Gamma))^{\prime}=L_{\infty,1/u}(\Gamma)$ and $(L_{\infty,u}(\Gamma))^{\prime}=L_{1,1/u}(\Gamma)$ with equality of norms.

Let $v,w\in L_0$ be weight functions on $\Gamma$ and $\Omega$, respectively, where $\mu(\Gamma)>0$. We consider the following spaces (see introduction),
\begin{equation*}
L_{1,w}(\Omega)\cap L_{\infty,v}(\Gamma)=\{x\in L_0:x\in L_{1,w}(\Omega)\text{ and }x\chi_{\Gamma}\in L_{\infty,v}(\Gamma)\}
\end{equation*}
equipped with the norm 
\begin{equation*}
\|x\|_{w,v}=\|x\|_{L_{1,w}(\Omega)\cap L_{\infty,v}(\Gamma)}=\max\{\|x\|_{1,w},\|x\chi_{\Gamma}\|_{\infty,v}\}\text{,}
\end{equation*}  
and
\begin{equation*}
L_{\infty,w}(\Omega)+L_{1,v}(\Gamma)=\{x\in L_0:x=y+z\text{ for some }y\in L_{\infty,w}(\Omega)\text{ and }z\in L_{1,v}(\Gamma)\}
\end{equation*}
equipped with the norm 
\begin{equation*}
\|x\|_{w,v}^{\Sigma}=\|x\|^{\Sigma}_{L_{\infty,w}(\Omega)+L_{1,v}(\Gamma)}=\inf\{\|y\|_{\infty,w}+\|z\|_{1,v}:x=y+z, y\in L_{\infty,w}(\Omega), z\in L_{1,v}(\Gamma)\} \text{.}
\end{equation*}  
Both $L_{1,w}(\Omega)\cap L_{\infty,v}(\Gamma)$ and $L_{\infty,w}(\Omega)+L_{1,v}(\Gamma)$  with their respective norms are Banach function lattices on $(\Omega,\Sigma,\mu)$ with the Fatou property and with $\supp(L_{1,w}(\Omega)\cap L_{\infty,v}(\Gamma))=\supp(L_{\infty,w}(\Omega)+L_{1,v}(\Gamma))=\Omega$. Moreover, the following K\"{o}the duality holds true (see Theorem \ref{dualXY}),
\begin{equation*}
(L_{1,w}(\Omega)\cap L_{\infty,v}(\Gamma))^{\prime}=L_{\infty,1/w}(\Omega)+L_{1,1/v}(\Gamma)\text{,}
\end{equation*}
and 
\begin{equation*}
(L_{\infty,1/w}(\Omega)+L_{1,1/v}(\Gamma))^{\prime}=L_{1,w}(\Omega)\cap L_{\infty,v}(\Gamma)\text{.}
\end{equation*}

The following lemma gives two useful conditions equivalent to the Daugavet property.  

\begin{Lemma}[See \cite{MR1621757} Lemma 2.2]
\label{Daug}
The following are equivalent.
\begin{enumerate}[{\rm(i)}]
\item A Banach space $(X,\|\cdot\|)$ has the Daugavet property.
\item\label{Daugii} For every $x\in S(X)$ and $y^*\in S(X^*)$ and every $\epsilon>0$ there is $x^*\in S(X^*)$ such that $x^*(x)>1-\epsilon$ and $\|x^*+y^*\|>2-\epsilon$.
\item\label{Daugiii} For every $x\in S(X)$ and $x^*\in S(X^*)$ and every $\epsilon>0$ there is $y\in S(X)$ such that $x^*(y)>1-\epsilon$ and $\|x+y\|>2-\epsilon$.
\end {enumerate}
\end{Lemma}

We need two technical lemmas.

\begin{Lemma}
\label{L1wpLivnorm}
Let $\Gamma$ be a measurable subset of $\Omega$ with $\mu(\Gamma)>0$ and $v,w\in L_0$ be weight functions on $\Gamma$ and $\Omega$, respectively. 
Let the space $X=L_{1,v}(\Gamma)+L_{\infty,w}(\Omega)$ be equipped with the norm $\|x\|=\|x\|_{w,v}^{\Sigma}$ 
Then 
\begin{equation*}
X^*\simeq (L_{\infty,1/v}(\Gamma)\cap L_{1,1/w}(\Omega))\oplus S
\end{equation*}
and the norm on $X^*$ is given by
\begin{equation*}
\|F\|=\max\{ \|f\chi_{\Gamma}\|_{\infty,1/v},\|f\|_{1,1/w}+\|T\| \}\text{,}
\end{equation*}
where $F=F_1+T$, $F_1(y)=\int_{\Omega}fy\,d\mu$ for some $f\in L_{\infty,1/v}(\Gamma)\cap L_{1,1/w}(\Omega)$ and $T$ is a singular functional on $X$.
\begin{proof}
The method of proof is similar to the proof of Theorem 2.12 \cite{MR3013764}. By Theorem \ref{dualXY} and from the general theory of Banach lattices we have that $X^*\simeq (L_{\infty,1/v}(\Gamma)\cap L_{1,1/w}(\Omega))\oplus S$. Let $x\in X$, that is $x=y+z$ where $y\in L_{1,v}(\Gamma)$ and $z\in L_{\infty,w}(\Omega)$. Let $F$, $F_1$, $T$ and $f$ be as in the statement of the theorem. Since $y$ is an order continuous element of $X$ we have that $T(y)=0$. Hence
\begin{equation*}
\begin{split}
|F(x)|&\leqslant |F_1(y)|+|F_1(z)|+|T(z)|\leqslant \int_{\Omega}|fy\chi_{\Gamma}|\,d\mu+\int_{\Omega}|fz|\,d\mu+\|T\|\|z\|_X\\
&\leqslant\|f\chi_{\Gamma}\|_{\infty,1/v}\|y\|_{1,v}+\|f\|_{1,1/w}\|z\|_{\infty,w}+\|T\|\|z\|_{\infty,w}\\
&\leqslant\max\{\|f\chi_{\Gamma}\|_{\infty,1/v},\|f\|_{1,1/w}+\|T\|  \}(\|y\|_{1,v}+\|z\|_{\infty,w})\text{.}
\end{split}
\end{equation*}
Therefore $\|F\|\leqslant\max\{\|f\chi_{\Gamma}\|_{\infty,1/v},\|f\|_{1,1/w}+\|T\|  \}$.

For every $0\leqslant x=y+z\in X$ with $y\in L_{1,v}(\Gamma)$ and $z\in L_{\infty,w}(\Omega)$, since $y$ is an order continuous element of $X$, $|T|(y)=0$. Hence, for every $\epsilon>0$, there are  $0\leqslant y_0\in L_{1,v}(\Gamma)$ and $0\leqslant z_0\in L_{\infty,w}(\Omega)$ such that
\begin{equation*}
\|y_0\|_{1,v}+\|z_0\|_{\infty,w}<1+\epsilon\text{ and } |T|(z_0)>\|T\|-\epsilon/2\text{.} 
\end{equation*} 
Since $f\in L_{1,1/w}(\Omega)$ there exists $0\leqslant z_1\in B(L_{\infty,w}(\Omega))$ such that 
\begin{equation*}
|F_1|(z_1)=\int_{\Omega}|f|z_1\, d\mu>\|f\|_{1,1/w}-\epsilon/2\text{.}
\end{equation*}
For $z_2=\max\{z_0,z_1\}$ we have that $0\leqslant z_2\in X$ and
\begin{equation*}
\begin{split}
|F|(z_2)&=|F_1|(z_2)+|T|(z_2)=\int_{\Omega}|f|z_2d\mu+|T|(z_2)\\
&\geqslant \int_{\Omega}|f|z_1\,d\mu+|T|(z_0)>\|f\|_{1,1/w}+\|T\|-\epsilon\text{.}
\end{split}
\end{equation*}
Since $\|z_2\|_{w,v}^{\Sigma}<1+\epsilon$, we get that
\begin{equation*}
\|F\|\geqslant(\|f\|_{1,1/w}+\|T\|-\epsilon)/(1+\epsilon)\text{.}
\end{equation*}
It follows that $\|F\|\geqslant\|f\|_{1,1/w}+\|T\|$.
We also have that
\begin{equation*}
\begin{split}
\|F\|&=\sup_{\|x\|_{w,v}^{\Sigma}\leqslant 1}|F(x)|\geqslant \sup_{\|x\chi_{\Gamma}\|_{1,v}\leqslant 1}|F(x)|\\
&=\sup_{\|x\chi_{\Gamma}\|_{1,v}\leqslant 1}|F_1(x)|=\sup_{\|x\chi_{\Gamma}\|_{1,v}\leqslant 1}\left|\int_{\Omega}f\chi_{\Gamma}x\,d\mu\right|=\|f\chi_{\Gamma}\|_{\infty,1/v}\text{.}
\end{split}
\end{equation*}
Hence $\|F\|\geqslant\max\{\|f\chi_{\Gamma}\|_{\infty,1/v},\|f\|_{1,1/w}+\|T\|  \} $ and the claim follows.
\end{proof}
\end{Lemma}

\begin{Lemma}
\label{L1wpLivOC}
 Let $\Gamma$ be a measurable subset of $\Omega$ with $\mu(\Gamma)>0$. Let $v,w\in L_0$ be weight functions on $\Gamma$ and $\Omega$, respectively. 
The space $X=L_{1,v}(\Gamma)+L_{\infty,w}(\Omega)$ equipped with the norm $\|x\|=\|x\|_{w,v}^{\Sigma}$ is order continuous if and only if $\mu(\Omega\setminus\Gamma)=0$ and $\int_{\Omega}v/w\, d\mu<\infty$.
\begin{proof}
Assume first that $\mu(\Omega\setminus\Gamma)=0$ and $\int_{\Omega}v/w\, d\mu<\infty$. 
 Let $x\in X$ be arbitrary, that is $x=y+z$ where $y\in L_{1,v}$ and $z\in L_{\infty,w}$. Since $\|z\|_{1,v}=\int_{\Omega}|z|wv/w\, d\mu\leqslant\|z\|_{\infty,w}\int_{\Omega}v/w\, d\mu<\infty$ we get $\|x\|_{1,v}<\infty$, that is $x\in L_{1,v}$. It follows that $X=L_{1,v}$ as sets and $\|x\|\leqslant\|x\|_{1,v}$ for every $x\in X$, and consequently $(X,\|\cdot\|)$ is order continuous. 

If $\mu(\Omega\setminus\Gamma)>0$ then in order to see that $X$ is not order continuous it is enough to take $x=(1/w)\chi_{\Omega\setminus\Gamma}$ and $x_n=(1/w)\chi_{A_n}$, $n\in\N$, where $(A_n)\subset\Omega\setminus\Gamma$ is a sequence of measurable sets such that $A_{n+1}\subset A_n$, $n\in\N$ and  $\mu(A_n)\to 0$ as $n\to\infty$. 

Assume now that $\mu(\Omega\setminus\Gamma)=0$ and  $\int_{\Omega}v/w\, d\mu=\infty$. We have two cases.  

Case 1. There is a measurable set $A\subset\Omega$ of finite measure such that $\int_{A}v/w\, d\mu=\infty$. In this case there exists a sequence $(A_n)$ of measurable subsets of $A$ such that $A_{n+1}\subset A_n$ for $n\in\N$, $\mu(A_n)\to 0$ as $n\to\infty$ and $\int_{A_n}v/w\, d\mu=\infty$, $n\in\N$. Let $x=(1/w)\chi_A$ and $x_n=(1/w)\chi_{A_n}$, $n\in\N$. We have that $0\leqslant x_n\leqslant x$ on $\Omega$ and $x_n\downarrow 0$ $\mu$-a.e. on $\Omega$. But $\|x\|=\|x_n\|=1$ for every $n\in\N$. Indeed, clearly $\|x_n\|\leqslant 1$. Let $x_n=y_n+z_n$ where $y_n\in L_{1,v}(\Gamma)$, $z_n\in L_{\infty,w}$ and for $c\in(0,1)$ denote $B_c=\{t\in A_n:|y_n(t)|>c/w(t)\}$. Since $y_n\in L_{1,v}(\Gamma)$ it must be that $\int_{B_c}v/w\, d\mu<\infty$. It follows that $\mu(A_n\setminus B_c)>0$, $|z_n|\geqslant(1-c)/w$ on $A_n\setminus B_c$ and $\|y_n\|_{1,v}+\|z_n\|_{\infty,w}\geqslant\|z_n\|_{\infty,w}\geqslant 1-c$. Since $c\in(0,1)$ was arbitrary, we get that $\|x_n\|=1$, $n\in\N$. Similarly $\|x\|=1$. Hence $X$ is not order continuous.

Case 2. For every measurable set $A\subset\Omega$ of finite measure $\int_Av/w\, d\mu<\infty$. In this case $\mu(\Omega)=\infty$. By $\sigma$-finiteness of $\mu$ there exists a sequence $(A_n)$ of measurable sets such that $\Omega=\cup_{n=1}^\infty A_n$, $A_n\subset A_{n+1}$, $\mu(A_n)<\infty$, $n\in\N$. Clearly $\int_{\Omega\setminus A_n}v/w\, d\mu=\infty$, $n\in\N$. Let $x=(1/w)\chi_{\Omega}$ and $x_n=(1/w)\chi_{\Omega\setminus A_n}$, $n\in\N$. We have that $0\leqslant x_n\leqslant x$ on $\Omega$ and $x_n\downarrow 0$ $\mu$-a.e. on $\Omega$. But $\|x\|=\|x_n\|=1$ for every $n\in\N$. Indeed, clearly $\|x_n\|\leqslant 1$. Let $x_n=y_n+z_n$ where $y_n\in L_{1,v}(\Gamma)$, $z_n\in L_{\infty,w}$ and for $c\in(0,1)$ denote $B_c=\{t\in \Omega\setminus A_n:|y_n(t)|>c/w(t)\}$. Since $y_n\in L_{1,v}(\Gamma)$ it must be that $\int_{B_c}v/w\, d\mu<\infty$. It follows that $\mu((\Omega\setminus A_n)\setminus B_c)>0$, $|z_n|\geqslant(1-c)/w$ on $(\Omega\setminus A_n)\setminus B_c$ and $\|y_n\|_{1,v}+\|z_n\|_{\infty,w}\geqslant\|z_n\|_{\infty,w}\geqslant1-c$. Thus $\|x_n\|=1$, $n\in\N$. Analogously we show that $\|x\|=1$. Hence $X$ is not order continuous.
\end{proof} 
\end{Lemma}

The following lemma gives conditions on weights $v$ and $w$ for which the space $L_{1,v}(\Gamma)+L_{\infty,w}(\Omega)$ fails the Daugavet property.

\begin{Lemma}
\label{L1wpLiv}
Let $\Gamma$ be a measurable subset of $\Omega$ with $\mu(\Gamma)>0$ and $v,w\in L_0$ be weight functions on $\Gamma$ and $\Omega$, respectively. 
The space $X=L_{1,v}(\Gamma)+L_{\infty,w}(\Omega)$ equipped with the norm $\|x\|=\|x\|_{w,v}^{\Sigma}$ does not have the Daugavet property whenever $\mu(\Omega\setminus\Gamma)>0$ or $\int_{\Gamma}v/w\, d\mu>1$.
\begin{proof}

We consider two cases. 

Case 1. Suppose first that $\mu(\Omega\setminus\Gamma)>0$ or $\int_{\Gamma}v/w\, d\mu=\infty$. 
By Lemma \ref{L1wpLivOC} the space $X$ is not order continuous. 

Clearly, there exist constants $\alpha,\beta>0$ and a measurable set $\Omega_0\subset \Gamma$ of positive and finite measure such that $\alpha\leqslant w(t),v(t)\leqslant\beta$ for all $t\in\Omega_0$. 
Let $A\subset\Omega_0$ be a measurable set such that $0<\mu(A)\leqslant\beta/\alpha$. It follows that 
\begin{equation}
\label{setA}
\frac{1}{\mu(A)}\frac{w(t)}{v(t)}\geqslant 1\text{ for all }t\in A\text{.}
\end{equation}
Define 
\begin{equation*}
x=\frac{1}{\mu(A)}\frac{1}{v}\chi_A\text{.}
\end{equation*}
Since $\|x\|_{1,v}=1$ we have that $\|x\|\leqslant 1$. Let $F_0\in X^*$ be defined by the function $f_0=v\chi_A$, that is $F_0(z)=\int_Avz\,d\mu$, $z\in X$. Since $f_0\in X^{\prime}=L_{\infty,1/v}(
\Gamma)\cap L_{1,1/w}(\Omega)$, $\|f_0\|_{\infty,1/v}=1$ and by (\ref{setA})
\begin{equation*}
\|f_0\|_{1,1/w}=\int_A\frac{v}{w}\,d\mu\leqslant 1\text{,}
\end{equation*}
we get that $\|F_0\|\leqslant 1$. Moreover $F_0(x)=1$. It follows that $\|x\|=1$.

Fix $b\in[1/2,1)$ and let $0<\epsilon<1$ and $c>1$ be such that
\begin{equation}
\label{ceps}
0<\frac{c\epsilon}{1-c\epsilon}\leqslant\mu(A)\alpha/\beta\quad\text{and}
\end{equation}  
\begin{equation}
\label{cepsb}
c\epsilon<1-b\text{.}
\end{equation}  

Define
\begin{equation*}
g=-bv\chi_{A}\text{.}
\end{equation*}
Clearly $\|g\|_{\infty,1/v}=b$ and by (\ref{setA}), $\|g\|_{1,1/w}\leqslant b$. Since the space $X$ is not order continuous, there are non-trivial singular functionals on $X$. Let $S_1\in B(X^*)$ be a singular functional on $X$ with the norm $\|S_1\|=1-\|g\|_{1,1/w}$. Define $G\in X^*$ by
\begin{equation*}
G(y)=\int_{\Omega}gy\, d\mu+S_1(y)\text{, }y\in X\text{.}
\end{equation*}
By Lemma \ref{L1wpLivnorm}, $\|G\|=\max\{\|g\|_{\infty,1/v},\|g\|_{1,1/w}+\|S_1\|\}=1$. 

Let $F\in S(X^*)$ be such that $F(x)>1-\epsilon$. Since $x\in L_{1,v}(\Gamma)$ and $L_{1,v}(\Gamma)$ is order continuous and $\|z\|\leqslant \|z\|_{1,v}$ for all $z\in L_{1,v}(\Gamma)$ it is clear that $x$ is an order continuous element of $X$. In view of $F=H+S_2$ where $H\in B(X^*)$ is an integral functional defined by a function $h\in B(L_{\infty,1/v}(\Gamma)\cap L_{1,1/w}(\Omega))$ and $S_2\in B(X^*)$ is a singular functional on $X$ and $x$ is an order continuous element of $X$, we get that $S_2(x)=0$. Hence
\begin{equation*}
1-\epsilon<F(x)=H(x)=\frac{1}{\mu(A)}\int_A\frac{h}{v}\, d\mu\text{.}
\end{equation*}

It is not difficult to see that for every $d>1$ there exists a  measurable subset $B\subset A$ with $\mu(B)>d^{-1}\mu(A)$ and
\begin{equation*}
h>(1-d\epsilon)v\text{ } \mu\text{-a.e. on }B\text{.}
\end{equation*}

 Taking $d=c$ and the corresponding $B\subset A$ from the above statement, in view of $v\geqslant\alpha$ and $w\leqslant\beta$ on $A$, we get that 
\begin{equation*} 
\|h\chi_B\|_{1,1/w}\geqslant (1-c\epsilon)\mu(B)\alpha/\beta\text{.}
\end{equation*}
Now, by Lemma \ref{L1wpLivnorm} we get $\|h\|_{1,1/w} + \|S_2\| \le \|F\| = 1$, and so 
\begin{equation*}
\|h\chi_{\Omega\setminus B}\|_{1,1/w}+\|S_2\|=\|h\|_{1,1/w}-\|h\chi_B\|_{1,1/w}+\|S_2\|\leqslant 1-(1-c\epsilon)\mu(B)\alpha/\beta\text{.}
\end{equation*} 

Now we show that $\|F+G\|_{X^{*}}\leqslant 2-\epsilon$. 
 Since $B\subset\Gamma$, $\|h\chi_{\Gamma}\|_{\infty,1/v}=\|(h/v)\chi_{\Gamma}\|_{\infty}\leqslant\|F\|=1$, $h/v>1-c\epsilon>b$ $\mu$-a.e. on $B$ by (\ref{cepsb}), and $1-b\leqslant b$ we get that
\begin{equation*}
\begin{split}
\|h+g\|_{1,1/w}+\|S_1+S_2\|&\leqslant\|(h+g)\chi_B\|_{1,1/w}+\|(h+g)\chi_{\Omega\setminus B}\|_{1,1/w}+\|S_1\|+\|S_2\|\\
&\leqslant\int_B\left|\frac{h}{v}-b\right|\frac{v}{w}\, d\mu+\|h\chi_{\Omega\setminus B}\|_{1,1/w}+\|S_2\|+\|g\chi_{\Omega\setminus B}\|_{1,1/w}+\|S_1\|  \\
&\leqslant (1-b)\left\|v\chi_B\right\|_{1,1/w}+1-(1-c\epsilon)\mu(B)\alpha/\beta + b\left\|v\chi_{A\setminus B}\right\|_{1,1/w}+\|S_1\|\\
&\leqslant1-(1-c\epsilon)\mu(B)\alpha/\beta+\|g\chi_A\|_{1,1/w}+\|S_1\|\\
&\leqslant2-(1-c\epsilon)\mu(B)\alpha/\beta\leqslant 2-\epsilon\text{,}
\end{split}
\end{equation*}
where the last inequality follows from $\mu(B)>c^{-1}\mu(A)$ and (\ref{ceps}). Moreover, by (\ref{cepsb}) we get that  
\begin{equation*}
\|(h+g)\chi_{\Gamma}\|_{\infty,1/v}\leqslant\|h\chi_{\Gamma}\|_{\infty,1/v}+\|g\|_{\infty,1/v}\leqslant 1+b\leqslant 2-\epsilon\text{.}
\end{equation*}
It follows that
\begin{equation*}
\|F+G\|_{X^{*}}=\max\{\|(h+g)\chi_{\Gamma}\|_{\infty,1/v},\|h+g\|_{1,1/w}+\|S_1+S_2\|\}\leqslant 2-\epsilon\text{,}
\end{equation*}
hence $X$ fails the Daugavet property by Lemma \ref{Daug}(\ref{Daugii}).

Case 2. Suppose now that $\mu(\Omega\setminus\Gamma)=0$ and $1<\int_{\Omega}v/w\, d\mu<\infty$. By Lemma \ref{L1wpLivOC} the space $X$ is order continuous. In this case the whole proof above can be repeated with the following modifications. If $\int_{\Omega}v/w\, d\mu\leqslant 2$ then fix $b=(\int_{\Omega}v/w\, d\mu)^{-1}$, $g=-bv\chi_{\Omega}$, and $S_1=S_2=0$. Otherwise, find a measurable subset $C\subset\Omega$ containing $A$ such that $\int_Cv/w\, d\mu=2$, fix $b=1/2$, $g=-bv\chi_{C}$ and $S_1=S_2=0$. In both cases we have $\|g\|_{1,1/w}=1$, $\|g\|_{\infty,1/v}=b<1$. Hence $\|G\|=\max\{\|g\|_{\infty,1/v},\|g\|_{1.1/w}\}=1$.
\end{proof}
\end{Lemma}

Now we can characterize $L_{1,v}(\Gamma)+L_{\infty,w}(\Omega)$ spaces with the Daugavet property.

\begin{Theorem}
\label{corPlus}
Let $\Gamma$ be a measurable subset of $\Omega$ such that $\mu(\Gamma)>0$ and $v,w\in L_0$ be weight functions on $\Gamma$ and $\Omega$, respectively. 
Let the space $X=L_{1,v}(\Gamma)+L_{\infty,w}(\Omega)$ be equipped with the norm $\|x\|=\|x\|_{w,v}^{\Sigma}$. The following conditions are equivalent.
\begin{enumerate}[{\rm(i)}]
\item $X$ has the Daugavet property.
\item $X=L_{1,v}$.
\item $\mu(\Omega\setminus\Gamma)=0$ and $\int_{\Omega}v/w\, d\mu\leqslant 1$.
\end{enumerate} 
\begin{proof}
 If $\int_{\Omega}v/w\, d\mu\leqslant 1$ and $\Gamma=\Omega$ up to a set of measure zero then $\|x\|_{1,v}=\int_{\Omega}|x|v\,d\mu=\int_{\Omega}|x|wv/w\,d\mu\leqslant\|x\|_{\infty,w}\int_{\Omega}v/w\,d\mu\leqslant\|x\|_{\infty,w}$ for every $x\in X$. It follows that, for every $x\in X$, if $x=y+z$ where $y\in L_{1,v}$ and $z\in L_{\infty,w}$ we have that $\|y\|_{1,v}+\|z\|_{\infty,w}\geqslant\|y\|_{1,v}+\|z\|_{1,v}\geqslant \|y+z\|_{1,v}=\|x\|_{1,v}$. Hence $\|x\|\geqslant\|x\|_{1,v}$. The opposite inequality is obvious.  Therefore $\|x\|=\|x\|_{1,v}$ for every $x\in X$. Since $L_{1,v}\subset X$ we have that $X=L_{1,v}$. Hence {\rm(iii)} implies {\rm(ii)} which in turn clearly implies  {\rm(i)}. By Lemma \ref{L1wpLiv}, {\rm(iii)} follows from {\rm(i)}.
\end{proof}
\end{Theorem}

Similarly as above we describe $L_{1,w}(\Omega)\cap L_{\infty,v}(\Gamma)$ spaces with the Daugavet property.
\begin{Lemma}
\label{L1wiLiv}
Let $\Gamma$ be a measurable subset of $\Omega$ such that $\mu(\Gamma)>0$ and $v,w\in L_0$ be weight functions on $\Gamma$ and $\Omega$, respectively. 
The space $X=L_{1,w}(\Omega)\cap L_{\infty,v}(\Gamma)$, with the norm $\|x\|=\|x\|_{w,v}$ does not have the Daugavet property whenever $\mu(\Omega\setminus\Gamma) > 0$ or $\int_{\Gamma}w/v\, d\mu>1$. 
\begin{proof}
We show that the condition (\ref{Daugiii}) of Lemma \ref{Daug} fails. We consider two cases.

{\rm(i)} Assume that $\mu(\Omega\setminus\Gamma) > 0$. 
 Let $c\in(0,1)$ and $A\subset\Gamma$ be a measurable set of positive and finite measure such that 
\begin{equation*}
c\int_A\frac{w}{v}\, d\mu=:\gamma\in(0,c)\text{.}
\end{equation*} 
Let $c_2>0$ and $A_2\subset\Omega\setminus\Gamma$ be a measurable set such that $c_2\int_{A_2}w\, d\mu=1-\gamma$. 
Define 
\begin{equation*}
x=c\frac{1}{v}\chi_{A}+c_2\chi_{A_2}\text{.}
\end{equation*}
Since $\|x\chi_{\Gamma}\|_{\infty,v}=c<1$ we have that $\|x\|=1=\|x\|_{1,w}$.

Let $F\in X^*$ be induced by $f=-(c/\gamma)w\chi_A$.
Since $X^{\prime}=L_{\infty,1/w}(\Omega)+L_{1,1/v}(\Gamma)$ and $\|f\|_{1,1/v}=1$ we get that $\|F\|\leqslant 1$. Let $z=(1/v)\chi_A$. Clearly $\|z\|_{\infty,v}=1$ and $\|z\|_{1,w}=\gamma/c<1$. Whence $\|z\|=1$. Moreover $F(z)=-1$. It follows that $\|F\|=1$. 

We finish the proof using methods similar to those in Theorem 2.9 \cite{MR3013764} or Proposition 4.3 \cite{MR2400252}.

Let $0<\epsilon<2\gamma(1-c)/(1-c+2\gamma)$ and $y\in S(X)$ be such that $F(y)>1-\epsilon$. Observe that $\epsilon<1-c$. Let 
\begin{equation*}
D=\{t\in A:-y(t)\geqslant 0\}\text{, }E=\left\{t\in D:-y(t)\leqslant \frac{c}{v(t)}\right\}\text{, }B=\left\{t\in D:-y(t)>\frac{c}{v(t)}\right\}\text{.}
\end{equation*}

Since $|yv|\leqslant 1$ $\mu$-a.e. on $\Gamma$, $D=E\cup B$ and $E\cap B=\emptyset$, we get that
\begin{equation*}
\begin{split}
1-\epsilon&<F(y)=\int_{A}\frac{c}{\gamma}w(-y)\, d\mu\leqslant \int_{D}\frac{c}{\gamma}w(-y)\, d\mu=\int_{E}\frac{c}{\gamma}w(-y)\, d\mu+\int_{B}\frac{c}{\gamma}w(-y)\, d\mu\\
&\leqslant\frac{c}{\gamma}\int_{E}c\frac{w}{v}\, d\mu+\frac{c}{\gamma}\int_{B}\frac{w}{v}\, d\mu\leqslant \frac{c}{\gamma}\left(c\int_{A}\frac{w}{v}\, d\mu-c\int_{B}\frac{w}{v}\, d\mu\right)+\frac{c}{\gamma}\int_{B}\frac{w}{v}\, d\mu\\
&=\frac{c}{\gamma}\left(\gamma-c\int_{B}\frac{w}{v}\, d\mu\right)+\frac{c}{\gamma}\int_{B}\frac{w}{v}\, d\mu=c+\frac{c}{\gamma}(1-c)\int_{B}\frac{w}{v}\,d\mu\text{.}
\end{split}
\end{equation*}
It follows that
\begin{equation}
\label{intF}
c\int_{B}\frac{w}{v}\, d\mu\geqslant \frac{\gamma}{1-c}(1-c-\epsilon)\text{.}
\end{equation}
We also have that
\begin{equation}
\label{xpyinf}
 \|(x+y)\chi_{\Gamma}\|_{\infty,v}\leqslant\|x\chi_{\Gamma}\|_{\infty,v}+\|y\chi_{\Gamma}\|_{\infty,v}\leqslant c+1<2-\epsilon\text{,}
\end{equation}
and
\begin{equation*}
\begin{split}
\|(x+y)\chi_D\|_{1,w}&=\|(x+y)\chi_E\|_{1,w}+\|(x+y)\chi_B\|_{1,w}\\
&=\int_E\left|\frac{c}{v}-(-y)\right|w\, d\mu+\int_B\left|\frac{c}{v}-(-y)\right|w\, d\mu\\
&=\int_E\left(\frac{c}{v}-(-y)\right)w\, d\mu+\int_B\left(-\frac{c}{v}-y)\right)w\, d\mu\\
&=c\left(\int_E\frac{w}{v}-\int_B\frac{w}{v}\right)+\int_B(-y)w\, d\mu-\int_E(-y)w\, d\mu\\ 
&=c\left(\int_E\frac{w}{v}-\int_B\frac{w}{v}\right)+\|y\chi_B\|_{1,w}-\|y\chi_E\|_{1,w}\text{.}
\end{split}
\end{equation*}
Moreover
\begin{equation*}
\begin{split}
\|x\chi_D\|_{1,w}+\|y\chi_D\|_{1,w}&=\|x\chi_E\|_{1,w}+\|x\chi_B\|_{1,w}+\|y\chi_E\|_{1,w}+\|y\chi_B\|_{1,w}\\
&=c\left(\int_E\frac{w}{v}\, d\mu+\int_B\frac{w}{v}\, d\mu\right)+\|y\chi_E\|_{1,w}+\|y\chi_B\|_{1,w}\text{.}
\end{split}
\end{equation*}
It follows that 
\begin{equation*}
\|x\chi_D\|_{1,w}+\|y\chi_D\|_{1,w}-\|(x+y)\chi_D\|_{1,w}=2c\int_{B}\frac{w}{v}\, d\mu+2\|y\chi_E\|_{1,w}\geqslant2c\int_{B}\frac{w}{v}\, d\mu\text{,}
\end{equation*}
hence
\begin{equation}
\label{xpyD}
\|(x+y)\chi_D\|_{1,w}\leqslant\|x\chi_D\|_{1,w}+\|y\chi_D\|_{1,w}-2c\int_{B}\frac{w}{v}\, d\mu\text{.}
\end{equation}
Finally by inequality (\ref{intF}) and definition of $\epsilon$, we get that
\begin{equation*}
\begin{split}
\|x+y\|_{1,w}&=\|(x+y)\chi_D\|_{1,w}+\|(x+y)\chi_{\Omega\setminus D}\|_{1,w}\\
&\leqslant\|x\chi_D\|_{1,w}+\|y\chi_D\|_{1,w}-2c\int_{B}\frac{w}{v}\, d\mu+\|x\chi_{\Omega\setminus D}\|_{1,w}+\|y\chi_{\Omega\setminus D}\|_{1,w}\\
&=2-2c\int_{B}\frac{w}{v}\, d\mu\leqslant 2-\frac{2\gamma}{1-c}(1-c-\epsilon)<2-\epsilon\text{.}
\end{split}
\end{equation*}
Whence by (\ref{xpyinf})
\begin{equation*}
\|x+y\|=\max\{\|x+y\|_{1,w},\|(x+y)\chi_{\Gamma}\|_{\infty,v}\}\leqslant 2-\epsilon\text{,}
\end{equation*}
which finishes the proof in this case by Lemma \ref{Daug}(\ref{Daugiii}).

{\rm(ii)} Assume now that $\mu(\Omega\setminus\Gamma)=0$ 
and $\int_{\Omega}w/v\, d\mu>1$. There are a set $A\subset\Omega$ with finite and positive measure and a constant $c\in(0,1)$ such that 
\begin{equation*}
c\int_{A}\frac{w}{v}\, d\mu=1\text{.}
\end{equation*}
Moreover, there is a measurable set $A_1\subset A$ such that 
\begin{equation*}
\int_{A_1}\frac{w}{v}\, d\mu=1\text{.}
\end{equation*}
Let
\begin{equation*}
x=c\frac{1}{v}\chi_{A}\text{,}
\end{equation*}
and $F\in X^*$ be induced by $f=-w\chi_{A_1}$. 
Since $X^{\prime}=L_{\infty,1/w}+L_{1,1/v}$ and $\|f\|_{\infty,1/w}=1$ we get that $\|F\|\leqslant 1$. Let $z=(1/v)\chi_{A_1}$. Clearly $\|z\|_{\infty,v}=\|z\|_{1,w}=1=\|z\|$. Moreover $F(z)=-1$. It follows that $\|F\|=1$. 

We finish the proof similarly as in the previous case. 

Let $0<\epsilon<\min\{2c(1-c)/(1+c),1-c\}$ and $y\in S(X)$ be such that $F(y)>1-\epsilon$. Let 
\begin{equation*}
D=\{t\in A_1:-y(t)\geqslant 0\}\text{, }E=\left\{t\in D:-y(t)\leqslant \frac{c}{v(t)}\right\}\text{, }B=\left\{t\in D:-y(t)>\frac{c}{v(t)}\right\}\text{.}
\end{equation*}

Since $|yv|\leqslant 1$ $\mu$-a.e. on $\Omega$, $D=E\cup B$ and $E\cap B=\emptyset$ we get that
\begin{equation*}
\begin{split}
1-\epsilon&<F(y)=\int_{A_1}w(-y)\, d\mu\leqslant \int_{D}w(-y)\, d\mu=\int_{E}w(-y)\, d\mu+\int_{B}w(-y)\, d\mu\\
&\leqslant\int_{E}c\frac{w}{v}\, d\mu+\int_{B}\frac{w}{v}\, d\mu\leqslant c\left(\int_{A_1}\frac{w}{v}\, d\mu-\int_{B}\frac{w}{v}\, d\mu\right)+\int_{B}\frac{w}{v}\, d\mu\\
&\leqslant c\left(1-\int_{B}\frac{w}{v}\, d\mu\right)+\int_{B}\frac{w}{v}\, d\mu=c+(1-c)\int_B\frac{w}{v}\, d\mu\text{.}
\end{split}
\end{equation*}
It follows that
\begin{equation}
\label{intB2}
c\int_{B}\frac{w}{v}\, d\mu\geqslant \frac{c}{1-c}(1-c-\epsilon)\text{.}
\end{equation}
In the same way as in case {\rm(i)} we obtain inequalities (\ref{xpyinf}) and (\ref{xpyD}).
Finally by inequalities (\ref{xpyD}), (\ref{intB2}) and definition of $\epsilon$, we get that
\begin{equation*}
\begin{split}
\|x+y\|_{1,w}&=\|(x+y)\chi_D\|_{1,w}+\|(x+y)\chi_{\Omega\setminus D}\|_{1,w}\\
&\leqslant\|x\chi_D\|_{1,w}+\|y\chi_D\|_{1,w}-2c\int_{B}\frac{w}{v}\, d\mu+\|x\chi_{\Omega\setminus D}\|_{1,w}+\|y\chi_{\Omega\setminus D}\|_{1,w}\\
&=2-2c\int_{B}\frac{w}{v}\, d\mu\leqslant 2-\frac{2c}{1-c}(1-c-\epsilon)<2-\epsilon\text{.}
\end{split}
\end{equation*}
Thus by (\ref{xpyinf})
\begin{equation*}
\|x+y\|=\max\{\|x+y\|_{1,w},\|x+y\|_{\infty,v}\}\leqslant 2-\epsilon\text{,}
\end{equation*}
which finishes the proof in this case by Lemma \ref{Daug}(\ref{Daugiii}).
\end{proof}
\end{Lemma}

\begin{Theorem}
\label{corInt}
Let $\Gamma$ be a measurable subset of $\Omega$ such that $\mu(\Gamma)>0$ and $v,w\in L_0$ be weight functions on $\Gamma$ and $\Omega$, respectively. 
Let $X=L_{1,w}(\Omega)\cap L_{\infty,v}(\Gamma)$ be equipped with the norm $\|x\|=\|x\|_{w,v}$. The following conditions are equivalent. 
\begin{enumerate}[{\rm(i)}]
\item $X$ has the Daugavet property.
\item $X=L_{\infty,v}$.
\item $\mu(\Omega\setminus\Gamma)=0$ and $\int_{\Omega}w/v\, d\mu\leqslant 1$.
\end{enumerate} 
\begin{proof}
Assuming {\rm(iii)} we have that $\|x\|_{1,w}=\int_{\Omega}|x|w\,d\mu=\int_{\Omega}|x|vw/v\,d\mu\leqslant\|x\|_{\infty,v}\int_{\Omega}w/v\, d\mu\leqslant\|x\|_{\infty,v}$. It follows that $\|x\|=\|x\|_{\infty,v}$ and therefore $X=L_{\infty,v}$. Hence {\rm(iii)} implies {\rm(ii)} which in turn clearly implies {\rm(i)}. By Lemma \ref{L1wiLiv} we have that {\rm(iii)} follows from {\rm(i)}.
\end{proof}
\end{Theorem}

\section{The Daugavet property in the Musielak-Orlicz spaces}

We begin this section with a basic observation regarding Orlicz functions.
\begin{Lemma}
\label{onehalf}
Let $\varphi$ be an Orlicz function. For every closed and bounded interval $I\subset(d_{\varphi},b_{\varphi})$ there is a constant $\sigma\in(0,1)$ such that $2\varphi(u/2)/\varphi(u)\leqslant \sigma$ for $u\in I$. Moreover, if $\varphi(b_{\varphi})<\infty$ then the same statement holds true for closed intervals $I\subset(d_{\varphi},b_{\varphi}]$.
\begin{proof}
It is not difficult to see that if $\varphi(u/2)=\varphi(u)/2$ for some $u\geqslant 0$ then $\varphi(v/2)=\varphi(v)/2$ for all $v\in[0,u]$.  Indeed, suppose that there are $u>0$ and $v\in(0,u)$ such that $\varphi(u/2)=\varphi(u)/2$ and $\varphi(v/2)<\varphi(v)/2$. Then $[\varphi(u/2)-\varphi(v/2)]/[u/2-v/2]>[\varphi(u)-\varphi(v)]/[u-v]$ which contradicts the convexity of $\varphi$. It follows that for all $u\in (d_{\varphi},b_{\varphi})$, $\varphi(u/2)<\varphi(u)/2$. 
Since the ratio $2\varphi(u/2)/\varphi(u)$ is a continuous function on $(d_{\varphi},b_{\varphi})$, for every closed and bounded interval $I\subset (d_{\varphi},b_{\varphi})$ there is a positive constant $\sigma<1$ such that $2\varphi(u/2)/\varphi(u)\leqslant \sigma$ for $u\in I$. In the case when $\varphi(b_{\varphi})<\infty$ the interval $I$ can include the point $b_{\varphi}$.
\end{proof}
\end{Lemma}

 The following lemma generalizes this fact to the case of a Musielak-Orlicz function. 

\begin{Lemma}
\label{L1}
Let $M$ be a Musielak-Orlicz function. If $\mu\{t\in \Omega:d_M(t)<b_M(t)\}>0$ then there exist a set $A\in\Sigma$ with positive and finite measure and numbers $a$, $b$ and $\sigma_1$ such that $a<b$, $0<\sigma_1<1$, $[a,b]\subset (d_M(t),b_M(t))$, $2M(t,u/2)/M(t,u)\leqslant\sigma_1$ for $t\in A$ and $u\in[a,b]$, and $M$ is bounded on $A\times[a,b]$.
\begin{proof}
Let $D=\{t\in \Omega:d_M(t)<b_M(t)\}$ and $\{r_1, r_2, \ldots\}\subset [0,\infty)$ be a countable dense set. Since $D=\cup_n\cup_m\{t\in D:d_M(t)<r_n<r_m<b_M(t)\}$, there are numbers $a<b$ and a set $D'\subset D$ such that $\mu D'>0$ and $[a,b]\subset(d_M(t),b_M(t))$ for $t\in D'$. 

Let $\sigma(t)=\sup\{2M(t,u/2)/M(t,u):u\in[a,b]\}$, $t\in D'$. 
Clearly $\sigma$ is a measurable function and $0<\sigma(t)<1$ for all $t\in D'$. Denote $D_n=\{t\in D':1-1/n\leqslant \sigma(t)<1-1/(n+1)\}$. Since $D'=\cup_nD_n$, $\mu(D_N)>0$ for some $N\in\N$. Define $\sigma_1=1-1/(N+1)$. Let $C_n=\{t\in D_N:n-1\leqslant M(t,b)<n\}$. Since $D_N=\cup_nC_n$ there is $N'\in\N$ such that $\mu C_{N'}>0$. Now, by taking as $A\in\Sigma$ any subset of $C_{N'}$ with positive and finite measure the claim follows.
\end{proof}
\end{Lemma}

Now we state a useful decomposition theorem for the Musielak-Orlicz spaces. The following notation will be used in the sequel. For a Musielak-Orlicz function $M$, we define
\begin{equation*}
\begin{split}
\Omega_{\infty}&=\{t\in \Omega:a_M(t)=b_M(t)\}\text{,}\quad v:\Omega\to[0,\infty)\text{, }\quad v=1/b_M\text{,}\\
\Omega_1&=\{t\in \Omega\setminus \Omega_{\infty}:0=a_M(t)<d_M(t)=b_M(t)=\infty\}\text{,}\\
\Omega_{1,\infty}&=\{t\in\Omega\setminus\Omega_{\infty}:0=a_M(t)<d_M(t)=b_M(t)<\infty\}\text{,}\\
&\quad w:\Omega\to[0,\infty)\text{,}\quad w(t)=\begin{cases}
M(t,u_t)/u_t\text{ where }u_t\in(0, b_M(t))\text{ if }t\in \Omega_1\cup\Omega_{1.\infty}\text{,}\\
0\text{ if }t\in \Omega\setminus (\Omega_1\cup\Omega_{1,\infty})\text{.}
\end{cases}
\end{split}
\end{equation*}
For $t\in\Omega_{\infty}$, $b_M(t)\in(0,\infty)$ and so $v(t)\in(0,\infty)$. Observe also that $w$ is well defined by the definitions of $\Omega_1$ and $\Omega_{1,\infty}$, and $w(t)\in(0,\infty)$ for $t\in\Omega_1\cup\Omega_{1,\infty}$. In fact $w=a_N$ $\mu$-a.e. on $\Omega$, where $N$ is the complementary function of $M$. Moreover, if $t\in \Omega_1\cup\Omega_{1,\infty}$ then $M(t,u)=w(t)u$ for all $u\in[0,b_M(t))$. Clearly $\mu(\Omega_{\infty}\cap (\Omega_1\cup\Omega_{1,\infty}))=0$.

It is easy to see that if $\Omega_{\infty}=\Omega$ up to a set of measure zero, then $L_M=L_{\infty,v}=L_{\infty,1/b_M}$ with $\|x\|_M=\|x\|_{\infty,v}$. If $d_M(t)=\infty$ for $\mu$-a.e. $t\in\Omega$ then $L_M=L_{1,w}=L_{1,a_N}$ with $\|x\|_M=\|x\|_{1,w}$. Moreover, if $\mu(\Omega_{\infty})=0$ and $d_M(t)=b_M(t)$ $\mu$-a.e. on $\Omega$ then 
\begin{equation*}
L_M=L_{\infty,v}(\Omega_{1,\infty})\cap L_{1,w}(\Omega)=L_{\infty,1/b_M}(\Omega_{1,\infty})\cap L_{1,a_N}(\Omega)
\end{equation*} and $\|x\|_M=\|x\|_{w,v}=\max\{\|x\|_{1,w},\|x\chi_{\Omega_{1,\infty}}\|_{\infty,v}\}$. 

\begin{Theorem}
\label{MOdecomp}
Let $M$ be a Musielak-Orlicz function. Then for any $x\in L_M$,
\begin{equation*}
\|x\|_M=\max\{\|x\chi_{\Omega_{\infty}}\|_{\infty,v},\|x\chi_{\Omega\setminus \Omega_{\infty}}\|_M\}\text{,}
\end{equation*}
 and thus
\begin{equation*}
L_M=L_{\infty,v}(\Omega_{\infty})\oplus_{\infty} L_M(\Omega\setminus \Omega_{\infty})\text{.}
\end{equation*}

Moreover, if $d_M(t)=b_M(t)$ $\mu$-a.e. on $\Omega\setminus \Omega_{\infty}$ then
\begin{equation*}
\|x\|_M=\max\{\|x\chi_{\Omega_{\infty}\cup\Omega_{1,\infty}}\|_{\infty,v},\|x\chi_{\Omega\setminus \Omega_{\infty}}\|_{1,w}\}\text{,}
\end{equation*}
and thus
\begin{equation*}
L_M=L_{\infty,v}(\Omega_{\infty})\oplus_{\infty} \left(L_{\infty,v}(\Omega_{1,\infty})\cap L_{1,w}(\Omega\setminus \Omega_{\infty})\right)\text{.}
\end{equation*}
\begin{proof}
It is easy to observe that, if $\lambda>0$ is such that $|x(t)|/\lambda>b_M(t)$ on a subset of $\Omega_{\infty}$ of positive measure, that is $\lambda<\|x\chi_{\Omega_{\infty}}\|_{\infty,v}$, then $\rho_M(x/\lambda)=\infty$. Moreover, if $\lambda>\|x\chi_{\Omega_{\infty}}\|_{\infty,v}$ then $\rho_M(x\chi_{\Omega_{\infty}}/\lambda)=0$. Hence
\begin{equation*}
\|x\|_M=\inf\{\lambda>\|x\chi_{\Omega_{\infty}}\|_{\infty,v}:\rho_M(x\chi_{\Omega\setminus \Omega_{\infty}}/\lambda)\leqslant 1\}=\max\{\|x\chi_{\Omega_{\infty}}\|_{\infty,v},\|x\chi_{\Omega\setminus \Omega_{\infty}}\|_M\}\text{.}
\end{equation*}

Assume now that $d_M(t)=b_M(t)$ $\mu$-a.e. on $\Omega\setminus \Omega_{\infty}$. 
Let $x\in L_M$ be such that $\mu(\Omega_{\infty}\cap\supp(x))=0$. 
 Similarly as above, if $\lambda<\|x\chi_{\Omega_{1,\infty}}\|_{\infty,v}$ then $\rho_M(x/\lambda)\geqslant\rho(x\chi_{\Omega_{1,\infty}}/\lambda)=\infty$. If $\lambda>\|x\chi_{\Omega_{1,\infty}}\|_{\infty,v}$ then $\rho_M(x/\lambda)=\rho_M(x\chi_{\Omega_1}/\lambda)+\rho(x\chi_{\Omega_{1,\infty}}/\lambda)=\int_{\Omega\setminus \Omega_{\infty}}w|x|/\lambda\, d\mu$. 
Hence
\begin{equation*}
\|x\|_M=\inf\left\{\lambda>\|x\chi_{\Omega_{1,\infty}}\|_{\infty,v}:\int_{\Omega\setminus \Omega_{\infty}}w|x|/\lambda\, d\mu\leqslant 1\right\}=\max\{\|x\chi_{\Omega_{1,\infty}}\|_{\infty,v},\|x\chi_{\Omega\setminus \Omega_{\infty}}\|_{1,w}\}\text{.}
\end{equation*}
For an arbitrary $x\in L_M$, applying the first part, we get that 
\begin{equation*}
\|x\|_M=\max\{\|x\chi_{\Omega_{\infty}}\|_{\infty,v},\|x\chi_{\Omega_{1,\infty}}\|_{\infty,v},\|x\chi_{\Omega\setminus \Omega_{\infty}}\|_{1,w}\}=\max\{\|x\chi_{\Omega_{\infty}\cup\Omega_{1,\infty}}\|_{\infty,v},\|x\chi_{\Omega\setminus \Omega_{\infty}}\|_{1,w}\}\text{.}
\end{equation*}
\end{proof}
\end{Theorem}

To characterize Musielak-Orlicz spaces with the Daugavet property we will use the following simple observation. A point $x\in S(X)$ is called \emph{uniformly non-$\ell_1^2$} (or \emph{uniformly non-square}), if there exists $\delta>0$ such that $\min(\|x+y\|,\|x-y\|)<2-\delta$ for all $y\in S(X)$. It is worth to mention that non-square points and non-squareness properties of the above type have been considered in context of many spaces \cite{MR2997393, MR2997764, MR3022822, MR1608225, MR2678885}.
 
\begin{Proposition}
\label{prop}
If $(X,\|\cdot\|)$ has the Daugavet property then there are no uniformly non-$\ell_1^2$ points either in $X$ or in $X^*$.
\begin{proof}
Let $x\in S(X)$ be arbitrary, $x^*\in S(X^*)$ be such that $x^*(x)=-1$ and $\epsilon>0$. By Lemma \ref{Daug}(iii), there is $y\in S(X)$ such that $x^*(y)>1-\epsilon$ and $\|x+y\|>2-\epsilon$. Also
$\|x-y\|=\|y-x\|\geqslant x^*(y-x)>2-\epsilon\text{.}$
Hence $x$ is not uniformly non-$\ell_1^2$.

Similarly, let $y^*\in S(X^*)$ and $\epsilon>0$ be arbitrary. There is $x\in S(X)$ such that $y^*(x)<-1+\epsilon/2$. By Lemma \ref{Daug}(ii), there is $x^*\in S(X^*)$ such that $x^*(x)>1-\epsilon/2$ and $\|x^*+y^*\|>2-\epsilon$. Also $\|x^*-y^*\|\geqslant (x^*-y^*)(x)>2-\epsilon$. Hence $y^*$ is not uniformly non-$\ell_1^2$.
\end{proof}
\end{Proposition}

The next proposition is known (Theorem 4 \cite{MR1602596}), but we present its proof for completeness. 
\begin{Proposition}
\label{LMLinf}
Let $M$ be a Musielak-Orlicz function. If $x\in L_M$ and $\rho_M(b_M\chi_{\supp(x)})\leqslant 1$ then $\|x\|_M=\|x\|_{\infty,v}$. Moreover, $L_M=L_{\infty,v}$ with equality of norms if and only if $\rho_M(b_M)\leqslant 1$.
\begin{proof}
It is clear that the condition $\rho_M(b_M\chi_{\supp(x)})\leqslant 1$ implies $b_M<\infty$ $\mu$-a.e. on $\supp(x)$. 

Suppose that $\|x\|_{\infty,v}=1$, that is $\|x/b_M\|_{\infty}=1$. Then for every $c>0$ there is a set of positive measure $D\subset\supp(x)$ such that $|x(t)|>(1-c)b_M(t)$, $t\in D$. It follows that $\rho_M(x/(1-c))=\infty$. Hence $\|x\|_M\geqslant 1$. Since $|x|\leqslant b_M$ $\mu$-a.e. on $\supp(x)$, $\rho_M(x)\leqslant\rho_M(b_M\chi_{\supp(x)})\leqslant 1$. It follows that $\|x\|_M=1$. 

Suppose now that $\|x\|_M=1$. It follows that $\rho_M(x)\leqslant 1$, $\rho_M(x/(1-c))>1$ for every $c>0$ and so $\|x\|_{\infty,v}\leqslant 1$. If $\|x\|_{\infty,v}=1-c$ for some $c>0$, that is $|x|\leqslant (1-c)b_M$ $\mu$-a.e. on $\Omega$, then $\rho_M(x/(1-c))\leqslant\rho_M(b_M\chi_{\supp(x)})\leqslant 1$, which gives a contradiction with $\|x\|_M=1$. Hence $\|x\|_{\infty,v}=1$.

We have showed that for all $x\in L_M$ with $\rho_M(b_M\chi_{\supp(x)})\leqslant 1$, $\|x\|_M=1$ if and only if $\|x\|_{\infty,v}=1$. The first claim follows. 

If $\rho_M(b_M)\leqslant 1$ then by the first part, it is clear that $L_M=L_{\infty,v}$ with equality of norms.
Let $L_M=L_{\infty,v}$ with equality of norms. Clearly $\|b_M\|_M=\|b_M\|_{\infty,v}=1$. It follows that  $\rho_M(b_M)\leqslant 1$.
\end{proof}
\end{Proposition}

\begin{Proposition}
\label{LMoL1}
Let $M$ and $N$ be complementary Musielak-Orlicz functions. If $x\in L_M^o$ and $\rho_N(b_N\chi_{\supp(x)})\leqslant 1$ then $\|x\|_M^o=\|x\|_{1,b_N}$. Moreover, $L_M^o=L_{1,b_N}$ with equality of norms if and only if $\rho_N(b_N)\leqslant 1$.
\begin{proof}
By the assumption $\rho_N(b_N\chi_{\supp(x)})\leqslant 1$ we have that $b_N<\infty$ $\mu$-a.e. on $\supp(x)$. Since for $\mu$-a.a. $t\in\Omega$, $\lim_{u\to\infty}M(t,u)/u=b_N(t)$ and $M(t,u)/u$ is increasing on $[0,\infty)$ \cite{MR1909821},
we have that $M(t,u)/u\leqslant b_N(t)$ for $\mu$-a.a. $t\in\Omega$. Hence $\|x\|_M^o=\inf_{k>0} k^{-1}[1+\rho_M(kx)]\leqslant\inf_{k>0} k^{-1}[1+\int_{\Omega}kb_N|x|\, d\mu]=\|x\|_{1,b_N}$.
On the other hand $\|x\|_M^o=\sup\{\int_{\Omega} xh\, d\mu:\rho_N(h)\leqslant 1\}\geqslant\int_{\Omega}xb_N\sign x\, d\mu=\|x\|_{1,b_N}$. The second statement follows trivially.
\end{proof}
\end{Proposition}

We also need the following result \cite[p. 64]{MR894273}.
\begin{Lemma}
\label{K_constr}
Let $M$ be a Musielak-Orlicz function with $b_M=\infty$ $\mu $-a.e. on $\Omega$. There exists an ascending sequence $(T_i)_{i=1}^{\infty}$ of measurable sets such that $\mu(T_i)<\infty$ for all $i\in\N$, $\mu(\Omega\setminus\cup_{i=1}^{\infty}T_i)=0$ and $\sup_{t\in T_i}M(t,u)<\infty$ for all $u\geqslant 0$ and all $i\in\N$.
\end{Lemma}

Before we proceed to the main theorem we need one more technical result. 
\begin{Lemma}
\label{LMa}
Let $M$ be a Musielak-Orlicz function. Then $E_M\ne\{0\}$ if and only if $\mu\{t\in \Omega:b_M(t)=\infty\}>0$.
\begin{proof}
It is clear that if $b_M<\infty$ $\mu$-a.e. on $\Omega$ then $E_M=\{0\}$. Now let $S=\{t\in \Omega:b_M(t)=\infty\}$. If $\mu(S)>0$, since $\mu$ is $\sigma$-finite, by Lemma \ref{K_constr}, $S=\cup_{i=1}^{\infty}S_i$, where $0<\mu(S_i)<\infty$, $i\in\N$ and for every $u\geqslant 0$, the function $M(\cdot,u)$ is essentially bounded on $S_i$. It follows that $\chi_{S_i}\in E_M$, $i\in\N$. 
\end{proof}
\end{Lemma}

The following theorem implies that a very wide class of Musielak-Orlicz spaces does not have the Daugavet property.

\begin{Theorem}
\label{nsqpt}
Let $M$ be a Musielak-Orlicz function. If $\mu\{t\in \Omega:d_M(t)<b_M(t)\}>0$ and $\rho_M(b_M)>1$ then there is a uniformly non-$\ell_1^2$ point in $L_M$. 
\begin{proof}
By Lemma \ref{L1} there exist a measurable set $C$ with $0<\mu C<\infty$ and numbers $a$, $b$ and $\sigma_1$ such that $a<b$, $0<\sigma_1<1$, $[a,b]\subset (d_M(t),b_M(t))$, 
\begin{equation}
\label{sigma}
M(t,u/2)\leqslant\sigma_1M(t,u)/2\text{ for }t\in C\text{ and }u\in[a,b]\text{,}
\end{equation}
 and $M$ is strictly positive and bounded on $C\times[a,b]$.

Without loss of generality we assume that $\rho_M(a\chi_C)\leqslant 1$. Denote $S=\{t\in \Omega:b_M(t)=\infty\}$.

We consider two cases.  First, suppose that $\mu(S)>0$. Then, let $A\subset C$ be such that $\mu(S\setminus A)>0$ and $\rho_M(a\chi_A)\leqslant 1$. By Lemma \ref{LMa}, $E_M(S\setminus A)\ne\{0\}$. By Lemma \ref{K_constr} there is a measurable set $T$ of positive and finite measure such that $\chi_T\in E_M(S\setminus A)$. Define $x=a\chi_A+x_0$, where $x_0=d_0\chi_T$ and $d_0\geqslant 0$ is such that  $\rho_M(a\chi_A+x_0)=1$.

Now, suppose that $\mu(S)=0$, that is $b_M<\infty$ $\mu$-a.e. on $\Omega$. 
Let $A\subset C$ be such that $\mu(A)>0$, $\rho_M(a\chi_A)\leqslant 1$ and $\rho_M(b_M\chi_{\Omega\setminus A})>1$. We can find a measurable set $G\subset \Omega\setminus A$ with positive and finite measure such that $\rho_M(b_M\chi_G)>1$. By the left continuity of the modular, we get that for some positive constant $c_1<1$, $\rho_M(c_1b_M\chi_G)>1$. 
Let $G_n=\{t\in G:n-1\leqslant M(t,c_1b_M(t))< n\}$, $n\in\N$. Since $\cup_{n=1}^{\infty}G_n=G$, for $r\in\N$ large enough, $\infty>\rho_M(c_1b_M\chi_{\cup_{n=1}^rG_n})>1$. Let $c_2>0$ be such that 
\begin{equation*}
\rho_M(x_0)=1-\rho_M(a\chi_A)\text{, where }x_0=\frac{c_1}{1+c_2}b_M\chi_{\cup_{n=1}^rG_n}\text{.}
\end{equation*}
Then $\rho_M(a\chi_A+x_0)=1$ and define again $x=a\chi_A+x_0$. 

In both cases $\rho_M(x)=1$ and so $\|x\|_M=1$. By the construction we have $\rho_M((1+\epsilon)a\chi_A)<\infty$ and $\rho_M((1+\epsilon)x_0)<\infty$ for some $\epsilon>0$. Since $a\chi_A$ and $x_0$ have disjoint supports,
\begin{equation}
\label{x1eps}
\rho_M((1+\epsilon)x)<\infty\text{ for some }\epsilon>0\text{.}
\end{equation}

Let $y\in S(L_M)$ be arbitrary. We split $A$ into disjoint union $A=A_1\cup A_2\cup A_3$, where $A_1=\{t\in A:b_M(t)=\infty\}$,  $A_2=\{t\in A:b_M(t)<\infty\text{ and }M(t,b_M(t))=\infty\}$ and $A_3=\{t\in A:M(t,b_M(t))<\infty\}$. 
Since for every $\lambda>1$, $\rho_M(y/\lambda)\leqslant 1$, there exist constants $c\in(0,1)$ and $d>0$ so large that $\mu(A\cap B)>0$, where
\begin{equation*}
B=\{t\in \Omega:|y(t)|\leqslant d\chi_{A_1}(t)+cb_M(t)\chi_{A_2}(t)+b_M(t)\chi_{A_3}(t)\}\text{.}
\end{equation*} 

Let $\sigma:\Omega\to[0,1)$ be defined by
\begin{equation*}
\sigma(t)=\begin{cases} \sup\left\{2M(t,u/2)/M(t,u):u\in [a,d] \right\}\text{, if }t\in A_1\\
\sup\left\{2M(t,u/2)/M(t,u):u\in [a,cb_M(t)] \right\}\text{, if }t\in A_2\\
\sup\left\{2M(t,u/2)/M(t,u):u\in [a,b_M(t)] \right\}\text{, if }t\in A_3\\
0\text{, otherwise.}
\end{cases}
\end{equation*}
It is not difficult to see that $\sigma$ is a finite measurable function and $\sigma>0$ on $A$. 
It follows that for every $t\in A$, 
\begin{equation*}
M(t,u/2)\leqslant\sigma(t)M(t,u)/2\text{,}
\end{equation*}  
 for all $u\in[a,d]$ if $t\in A_1$, for all $u\in[a,cb_M(t)]$ if $t\in A_2$, and for all $u\in[a,b_M(t)]$ if $t\in A_3$. Since 
\begin{equation*}
A\cap B=\cup_{n=1}^{\infty}\{t\in A\cap B:1-1/n<\sigma(t)\leqslant1-1/(n+1)\}\text{,}
\end{equation*}
there is a subset $H\subset A\cap B$, with $\mu(H)>0$ and a constant $\sigma_2\in(0,1)$ such that
\begin{equation}
\label{sigma_t}
M(t,u/2)\leqslant\sigma_2M(t,u)/2\text{,}
\end{equation}  
 for all $u\in[a,d]$ if $t\in H\cap A_1$, for all $u\in[a,cb_M(t)]$ if $t\in H\cap A_2$, and for all $u\in[a,b_M(t)]$ if $t\in H\cap A_3$.

Let $\sigma_0=\max\{\sigma_1,\sigma_2\}$, $\eta=\rho_M(a\chi_A)$ and $\gamma=\rho_M(a\chi_{A\setminus H})$. Clearly $\sigma_0\in(0,1)$, $\eta\in(0,1]$ and $\gamma\in[0,\eta)$. It follows that
\begin{equation}
\label{rhoaH}
\rho_M(a\chi_H)=\rho_M(a\chi_A)-\rho_M(a\chi_{A\setminus H})=\eta-\gamma>0\text{.}
\end{equation}
Let $\delta\in(0,(1-\sigma_0)(\eta-\gamma)/2)$, that is 
\begin{equation}
\label{sigmainq}
(1-\sigma_0)(\eta-\gamma)/2-\delta>0\text{.}
\end{equation}

Since $a<b$, by (\ref{x1eps}) we find $\epsilon>0$  such that $(1+\epsilon)a\leqslant b$,
$
\rho_M((1+\epsilon)x)<\infty\text{,}
$
and
\begin{equation}
\label{epsdel}
\rho_M((1+\epsilon)x)<\rho_M(x)+\delta\text{.}
\end{equation} 
Define $z=(1+\epsilon)x$.

We finish the proof in a similar way as in the proof of Theorem 3.27 \cite[p. 133]{MR1410390}.

Define $D=\{t\in H:x(t)y(t)\geqslant 0\}$, $E=H\setminus D$. 
Since $D\subset H\subset A\cap B$, by (\ref{sigma_t}) and by definition of the set $B$, we have that
\begin{equation*}
\rho_M((y/2)\chi_D)\leqslant(\sigma_2/2)\rho_M(y\chi_D)\text{.}
\end{equation*}
Moreover, since $D\subset A$, $z=(1+\epsilon)(a\chi_A+x_0)$ and $\supp(x_0)\cap A=\emptyset$ we have that $|z|\chi_D=(1+\epsilon)a\chi_A$ and hence by (\ref{sigma}), taking into account that $A\subset C$, we get  
\begin{equation*}
\rho_M((z/2)\chi_D)\leqslant(\sigma_1/2)\rho_M(z\chi_D)\text{.}
\end{equation*}
It follows that
\begin{equation*}
\rho_M\left(\frac{z-y}{2}\chi_D\right)\leqslant\rho_M\left(\frac{\max\{|z|,|y|\}}{2}\chi_D\right)\leqslant\frac{\sigma_0}{2}\rho_M(\max\{|z|,|y|\}\chi_D)\leqslant\frac{\sigma_0}{2}\left(\rho_M(z\chi_D)+\rho_M(y\chi_D)\right)\text{.}
\end{equation*}
Similarly
\begin{equation*}
\rho_M\left(\frac{z+y}{2}\chi_E\right)\leqslant\frac{\sigma_0}{2}\left(\rho_M(z\chi_E)+\rho_M(y\chi_E)\right)\text{.}
\end{equation*}
Hence
\begin{equation}
\label{rhozyH}
\begin{split}
\rho_M\left(\frac{z+y}{2}\chi_H\right)&+\rho_M\left(\frac{z-y}{2}\chi_H\right)=\rho_M\left(\frac{z+y}{2}\chi_D\right)+\rho_M\left(\frac{z+y}{2}\chi_E\right)\\
&+\rho_M\left(\frac{z-y}{2}\chi_D\right)+\rho_M\left(\frac{z-y}{2}\chi_E\right)\leqslant\frac{1}{2}\rho_M(z\chi_D)+\frac{1}{2}\rho_M(y\chi_D) \\
&+\frac{\sigma_0}{2}\left(\rho_M(z\chi_E)+\rho_M(y\chi_E)\right)+\frac{\sigma_0}{2}\left(\rho_M(z\chi_D)+\rho_M(y\chi_D)\right)\\
&+\frac{1}{2}\rho_M(z\chi_E)+\frac{1}{2}\rho_M(y\chi_E)=\frac{1+\sigma_0}{2}\left(\rho_M(z\chi_H)+\rho_M(y\chi_H)\right)\text{.}
\end{split}
\end{equation}
By the inequality (\ref{epsdel}) and definition of $z$,
\begin{equation*}
2+\delta\geqslant\rho_M(y)+\rho_M(x)+\delta\geqslant\rho_M(z)+\rho_M(y)\text{.}
\end{equation*}
By the above, (\ref{rhoaH}), (\ref{rhozyH}) and convexity of the functions $M(t,\cdot)$ for $\mu$-a.e. $t\in \Omega$, we conclude that 
\begin{equation*}
\begin{split}
&2+\delta-\rho_M\left(\frac{z+y}{2}\right)-\rho_M\left(\frac{z-y}{2}\right) \\
&\geqslant\rho_M(z)+\rho_M(y)-\rho_M\left(\frac{z+y}{2}\right)-\rho_M\left(\frac{z-y}{2}\right)\\
&\geqslant\rho_M(z\chi_H)+\rho_M(y\chi_H)-\rho_M\left(\frac{z+y}{2}\chi_H\right)-\rho_M\left(\frac{z-y}{2}\chi_H\right)\\
&\geqslant\frac{1-\sigma_0}{2}\left(\rho_M(z\chi_H)+\rho_M(y\chi_H)\right)\geqslant\frac{1-\sigma_0}{2}\rho_M(z\chi_H)\geqslant\frac{1-\sigma_0}{2}\rho_M(a\chi_H)=\frac{1-\sigma_0}{2}\left(\eta-\gamma\right)\text{.}
\end{split}
\end{equation*} 
By (\ref{sigmainq}) we get that
\begin{equation*}
2-\rho_M\left(\frac{z+y}{2}\right)-\rho_M\left(\frac{z-y}{2}\right)> 0\text{,}
\end{equation*}
thus
\begin{equation*}
\min\left\{\rho_M\left(\frac{z+y}{2}\right),\rho_M\left(\frac{z-y}{2}\right)\right\}\leqslant 1\text{.}
\end{equation*}
If $\rho_M\left(\frac{z+y}{2}\right)\leqslant 1$ then $\left\|\frac{z+y}{2}\right\|_M\leqslant 1$ which gives $\left\|\frac{x+y/(1+\epsilon)}{2}\right\|_M\leqslant\frac{1}{1+\epsilon}$. Hence
\begin{equation*}
\begin{split}
&\left|\left\|\frac{x+y}{2}\right\|_M-\left\|\frac{x+y/(1+\epsilon)}{2}\right\|_M\right|\leqslant\left\|\frac{x+y}{2}-\frac{x+y/(1+\epsilon)}{2}\right\|_M\\\
&=\left\|\frac{y}{2}-\frac{y}{2(1+\epsilon)}\right\|_M=\frac{\epsilon}{2(1+\epsilon)}\text{.}
\end{split}
\end{equation*}
Therefore
\begin{equation*}
\left\|\frac{x+y}{2}\right\|_M\leqslant\frac{1}{1+\epsilon}+\frac{\epsilon}{2(1+\epsilon)}=1-\frac{\epsilon}{2(1+\epsilon)}\text{.}
\end{equation*}
If $\rho_M\left(\frac{z-y}{2}\right)\leqslant 1$ then we get similarly that
\begin{equation*}
\left\|\frac{x-y}{2}\right\|_M\leqslant1-\frac{\epsilon}{2(1+\epsilon)}\text{.}
\end{equation*}
Finally, for all $y\in S(L_M)$,
\begin{equation*}
\min\{\|x+y\|_M,\|x-y\|_M\}\leqslant 2-\epsilon/(1+\epsilon)\text{,}
\end{equation*}
that is $x$ is a uniformly non-$\ell^2_1$ point.
\end{proof}
\end{Theorem}

From the above theorem and from Proposition \ref{prop} we have the following corollary.
\begin{Corollary}
\label{DCor}
Let $M$ be a Musielak-Orlicz function. If $\mu\{t\in \Omega:d_M(t)<b_M(t)\}>0$ and $\rho_M(b_M)>1$ then  $L_M$ does not have the Daugavet property.
\end{Corollary}

We need two more results before we state the main theorem.  The following lemma is analogous to Lemma 4.2 \cite{MR3013764} proved there for the maximum norm.  
\begin{Lemma}
\label{sum1daug}
Let $X=X_1\oplus_1 X_2\oplus_1\ldots\oplus_1 X_n$ be a finite direct sum of Banach spaces $(X_i,\|\cdot\|_{i})$, $i=1, 2, \ldots n$, equipped with the norm $\|x\|=\|x_1\|_1+\|x_2\|_2+\ldots+\|x_n\|_n$, where $x=(x_1, x_2, \ldots, x_n)$, $x_i\in X_i$, $i=1, 2, \ldots, n$. If $X$ has the Daugavet property then it is inherited by each component $X_i$, $i=1, 2, \ldots, n$.
\begin{proof}
Without loss of generality we assume that $n=2$. Suppose that $X=X_1\oplus_1 X_2$ has the Daugavet property. It is enough to show that $X_1$ has that property. Let $T(x_1)=x_1^*(x_1)y_1$, $x_1\in X_1$, be an arbitrary rank $1$ operator on $X_1$, where $x_1^*\in X_1^*$, $y_1\in X_1$ and $\|x_1^*\|=\|y_1\|_1=1$. Clearly $\|T\|=1$ and $\|I+T\|_{X_1\to X_1}\leqslant 2$. We will show the opposite inequality. For any $x=(x_1,x_2)\in X$ define 
\begin{equation*}
x^*(x)=x_1^*(x_1)\text{,}\quad y=(y_1,0)\text{ and }\tilde{T}(x)=x^*(x)y=(x_1^*(x_1)y_1,0)\text{.}
\end{equation*}
Since $X^*\simeq(X_1^*\oplus_{\infty} X_2^*)$, $\|x^*\|=\|x_1^*\|=1$. Moreover $\|y\|=\|y_1\|_1=1$, hence $\|\tilde{T}\|_{X\to X}=1$. Since $\tilde{T}$ is a rank one operator on $X$, by the Daugavet property of $X$, \begin{equation*}
\begin{split}
2&=\|I+\tilde{T}\|_{X\to X}=\sup_{\|x\|\leqslant 1}\|x+x^*(x)y\|\\
&=\sup_{\|x_1\|_1+\|x_2\|_2\leqslant 1}\left\{\|x_1+x_1^*(x_1)y_1\|_1+\|x_2\|_2\right\}\text{.}
\end{split}
\end{equation*}
Hence for every $\epsilon>0$ there is $x\in X$, $x=(x_1,x_2)$, $\|x\|=\|x_1\|_1+\|x_2\|_2\leqslant 1$ such that
\begin{equation}
\label{inq2}
\|x_1+x_1^*(x_1)y_1\|_1> 2-\epsilon-\|x_2\|_2\geqslant 1+\|x_1\|_1-\epsilon\text{.}
\end{equation}
It follows that 
\begin{equation*}
2\geqslant\left\|\frac{x_1}{\|x_1\|_1}+x_1^*\left(\frac{x_1}{\|x_1\|_1}\right)y_1\right\|_1>\frac{1}{\|x_1\|_1}+1-\frac{\epsilon}{\|x_1\|_1}\text{.}
\end{equation*}
Multiplying by $\|x_1\|_1$, we get that $\|x_1\|_1\geqslant 1-\epsilon$. Hence by (\ref{inq2}),
\begin{equation*}
\|I+T\|\geqslant\|x_1+T(x_1)\|_1=\|x_1+x_1^*(x_1)y_1\|_1>2-2\epsilon\text{.}
\end{equation*}
It follows that $\|I+T\|=2$. 
\end{proof} 
\end{Lemma}

Recall that a weak$^*$ slice of $B(X^*)$ is a set of the form $\{f\in B(X^*):f(x)>1-\epsilon\}$, where $x\in S(X)$ and $0<\epsilon<1$. Observe that any weak$^*$ slice of $B(X^*)$ is a slice of $B(X^*)$. For a proof of the following proposition see Proposition I.1.11 \cite{MR1211634} and Lemma 3.1 \cite{DualDiamJConvA}.
\begin{Proposition}
\label{2rough}
Let $X$ be a Banach space. The following conditions are equivalent.
\begin{enumerate}[{\rm(i)}]
\item The norm on $X$ is $2$-rough, that is for all $x\in X$,
\begin{equation*}
\limsup_{\|h\|\to 0}\frac{\|x+h\|+\|x-h\|-2\|x\|}{\|h\|}=2\text{.}
\end{equation*}
\item $X^*$ has the weak$^*$ slice (or local) diameter $2$ property, that is every weak$^*$ slice of $B(X^*)$ has diameter $2$.
\item The norm on $X$ is locally octahedral, that is no point of $S(X)$ is uniformly non-$\ell_1^2$. 
\end{enumerate}
\end{Proposition}

Finally, we state theorem which characterizes Musielak-Orlicz spaces with the Daugavet property.

\begin{Theorem}
\label{MODaug}
Let $M$ and $N$ be complementary Musielak-Orlicz functions, $v=1/b_M$, $w=a_N$, $\Omega_{\infty}=\{t\in\Omega:a_M(t)=b_M(t)\}$, $\Omega_1=\{t\in \Omega\setminus \Omega_{\infty}:d_M(t)=b_M(t)=\infty\}$ and $\Omega_{1,\infty}=\{t\in\Omega\setminus\Omega_{\infty}:d_M(t)=b_M(t)<\infty\}$.
 The following conditions are equivalent.
\begin{enumerate}[{\rm(i)}]
\item $L_M$ has the Daugavet property.
\item $L_M=L_{1,w}$ or $L_M=L_{\infty,v}$ or $L_M=L_{\infty,v}(\Omega_{\infty})\oplus_{\infty}L_{1,w}(\Omega\setminus \Omega_{\infty})$.
\item $L_N^o=L_{\infty,1/w}$ or $L_N^o=L_{1,1/v}$ or $L_N^o=L_{1,1/v}(\Omega_{\infty})\oplus_1 L_{\infty,1/w}(\Omega\setminus \Omega_{\infty})$.
\item $L_N^o$ has the Daugavet property.
\end{enumerate}
\begin{proof}
Let $L_M$ have the Daugavet property. By Corollary \ref{DCor}, we have that $\mu\{t\in \Omega:d_M(t)<b_M(t)\}=0$ or $\rho_M(b_M)\leqslant 1$. The latter condition is equivalent to $L_M=L_{\infty,v}$ by Proposition \ref{LMLinf}. The former condition implies that for $\mu$-a.e. $t\in \Omega\setminus\Omega_{\infty}$ we have $d_M(t)=b_M(t)$. 
 Hence in view of Theorem \ref{MOdecomp}, 
\begin{equation}
\label{LMcases}
L_M=\begin{cases}
L_{1,w}\text{, if }\mu(\Omega\setminus\Omega_1)=0\text{,}\\
L_{\infty,v}\text{, if }\mu(\Omega\setminus\Omega_{\infty})=0\text{,}\\
L_{\infty,v}(\Omega_{\infty})\oplus_{\infty}L_{1,w}(\Omega\setminus\Omega_{\infty})\text{, if }\mu(\Omega_{\infty}),\mu(\Omega\setminus\Omega_{\infty})>0\text{ and }\mu(\Omega_{1,\infty})=0\text{,}\\
L_{\infty,v}(\Omega_{\infty})\oplus_{\infty}(L_{1,w}(\Omega\setminus\Omega_{\infty})\cap L_{\infty,v}(\Omega_{1,\infty}))\text{, if }\mu(\Omega\setminus\Omega_{\infty}), \mu(\Omega_{1,\infty})>0\text{.}
\end{cases}
\end{equation}
By Theorem \ref{corInt} applied to $\Omega \setminus \Omega_\infty$ and $\Omega_{1,\infty}$ for $\Omega$ and $\Gamma$ respectively, the second component of the last space in (\ref{LMcases}) has the Daugavet property if and only if  it is equal to $L_{\infty,v}(\Omega\setminus\Omega_{\infty})$. Hence, in view of Lemma 4.2 \cite{MR3013764} we see that (i) implies (ii).
Since the Daugavet property is lifted from components of $\oplus_{\infty}$ sums to the whole space \cite{MR1126202} we conclude that conditions (i) and (ii) are equivalent. Since the last statement is also true for $\oplus_1$ sums \cite{MR1126202}, we see that (iii) implies (iv).
  
Assume now that (iv) holds true. We will show that $\mu\{t\in \Omega:d_M(t)<b_M(t)\}=0$ or $\rho_M(b_M)\leqslant 1$. Suppose that this condition is not satisfied. Then by Theorem \ref{nsqpt} we have that $L_M$ is not locally octahedral.
 Hence, by Proposition \ref{2rough} the dual space $(L_M)^*\simeq L_N^o\oplus S$ fails the weak$^*$ slice diameter $2$ property. Therefore we can find  a weak$^*$ slice 
\begin{equation*}
S(x,\epsilon)=\{f\in B((L_M)^*):f(x)>1-\epsilon\}
\end{equation*}
with the diameter less than $2$, where $x\in S(L_M)$ and $\epsilon>0$. Let $\kappa:L_M\to(L_M)^{**}$ be the canonical mapping defined by $(\kappa(x))(x^*)=x^*(x)$, $x^*\in X^*$. Consider the sets 
\begin{equation*}
S^{\prime}(x,\epsilon)=\{f\in B((L_M)^*):f\in (L_M)^{*}_c\text{ and }f(x)>1-\epsilon \}
\end{equation*}
and
\begin{equation*}
S^{\prime\prime}(F,\epsilon)=\{y\in B(L_N^o):F(y)>1-\epsilon\}
\end{equation*}
where $F=\kappa(x)$. Since $(L_M)^{*}_c\simeq (L_M)^{\prime}=L_N^o$, there is a bijective correspondence preserving norm between $S^{\prime}(x,\epsilon)$ and $S^{\prime\prime}(F,\epsilon)$. 
Since $S^{\prime}(x,\epsilon)\subset  S(x,\epsilon)$ we see that the slice $S^{\prime\prime}(F,\epsilon)$ of $B(L_N^o)$ has the diameter less than $2$. Hence $L_N^o$ fails the slice diameter $2$ property. In particular $L_N^o$ does not have the Daugavet property, which contradicts (iv).

Hence, indeed it must be that $\mu\{t\in \Omega:d_M(t)<b_M(t)\}=0$ or $\rho_M(b_M)\leqslant 1$. The latter condition is equivalent to $L_N^o=L_{1,1/v}$ by Proposition \ref{LMoL1}. Similarly as previously, the former condition together with Theorem \ref{MOdecomp}, the K\"{o}the duality $(L_M)^{\prime}=L_N^o$ (Theorem \ref{dualLM}) and (\ref{LMcases}) gives $L_N^o=L_{\infty,1/w}$, or $L_N^o=L_{1,1/v}$, or $L_N^o=L_{1,1/v}(\Omega_{\infty})\oplus_1 L_{\infty,1/w}(\Omega\setminus \Omega_{\infty})$, or
\begin{equation*}
L_N^o=L_{1,1/v}(\Omega_{\infty})\oplus_1(L_{1,1/v}(\Omega_{1,\infty})+L_{\infty,1/w}(\Omega\setminus \Omega_{\infty}))\text{,}
\end{equation*}
 where the norm on the second component is $\|x\|^{\Sigma}_{1/w,1/v}$. Since the latter space has the Daugavet property, by Lemma \ref{sum1daug}  we infer that the second component of that space has the Daugavet property as well. 
Now we see that the condition (iii) follows from Theorem \ref{corPlus}. Hence (iii) and (iv) are equivalent. 

Conditions (ii) and (iii) are clearly equivalent by the K\"othe duality $L_N^o=(L_M)^{\prime}$ and $(L_N^o)^{\prime}=L_M$ (see the appendix).
\end{proof}
\end{Theorem}

\begin{Corollary}
Let $M$ be a Musielak-Orlicz function such that $0<M(t,u)<\infty$ for $\mu$-a.a. $t\in\Omega$ and for all $u>0$, that is $a_M=0$ and $b_M=\infty$ $\mu$-a.e. on $\Omega$. Let $N$ be the function complementary to $M$.   
The following conditions are equivalent.
\begin{enumerate}[{\rm(i)}]
\item $L_M$ has the Daugavet property.
\item $L_M=L_{1,a_N}$. 
\item $L_N^o=L_{\infty,1/a_N}$.
\item $L_N^o$ has the Daugavet property.
\end{enumerate}
\end{Corollary}

As we noted in the introduction, if $M(t,u)=\varphi(u)$ for all $t\in\Omega$ and $u\geqslant 0$, where $\varphi$ is an Orlicz function then $L_M=L_{\varphi}$, the Orlicz space generated by $\varphi$. In this case $a_M=a_{\varphi}$ and $b_M=b_{\varphi}$ on $\Omega$, where $a_{\varphi}$ and $b_{\varphi}$ are constants defined in the introduction.  

\begin{Corollary}
\label{DaugOrl}
Let $\varphi$ and $\psi$ be complementary Orlicz functions. The following statements are equivalent.
\begin{enumerate}[{\rm(i)}]
\item $L_{\varphi}$ has the Daugavet property.
\item $L_{\varphi}=L_{1,a_{\psi}}$ or $L_{\varphi}=L_{\infty,1/b_{\varphi}}$.
\item $L_{\psi}^o=L_{\infty,1/a_{\psi}}$ or $L_{\psi}^o=L_{1,b_{\varphi}}$.
\item $L_{\psi}^o$ has the Daugavet property.
\end{enumerate}
\end{Corollary}

Another corollary from Theorem \ref{MODaug} is the following generalization of Theorem 4.1 \cite{MR3013764}.
\begin{Corollary}
Let $L_{p(t)}$ be a Nakano space, where $1\leqslant p(t)\leqslant\infty$ and  $1/p(t)+1/q(t)=1$ for $\mu$-a.a. $t\in \Omega$ with the usual convention that $q(t)=\infty$ if $p(t)=1$. Denote $\Omega_{\infty}=\{t\in\Omega:p(t)=\infty\}$. The following statements are equivalent.
\begin{enumerate}[{\rm(i)}]
\item $L_{p(t)}$ has the Daugavet property.
\item $L_{p(t)}=L_{1}$ or $L_{p(t)}=L_{\infty}$ or $L_{p(t)}=L_{1}(\Omega\setminus\Omega_{\infty})\oplus_{\infty}L_{\infty}(\Omega_{\infty})$.
\item $L_{q(t)}^o=L_{\infty}$ or $L_{q(t)}^o=L_{1}$ or $L_{q(t)}^o=L_{\infty}(\Omega\setminus\Omega_{\infty})\oplus_1 L_{1}(\Omega_{\infty})$.
\item $L_{q(t)}^o$ has the Daugavet property.
\end{enumerate}
\end{Corollary}

From the proof of Theorem \ref{MODaug} we can also deduce the following result.
\begin{Corollary}
If $L_N^o$ has the slice diameter $2$ property then $L_N^o=L_{\infty,1/a_N}$, or $L_N^o=L_{1,b_M}$, or $L_N^o=L_{1,b_M}(\Omega_{\infty})\oplus_1 L_{\infty,1/a_N}(\Omega\setminus \Omega_{\infty})$, or $L_N^o=L_{1,b_M}(\Omega_{\infty})\oplus_1 (L_{1,b_M}(\Omega_{1,\infty})+L_{\infty,1/a_N}(\Omega\setminus \Omega_{\infty}))$.
\end{Corollary}

\section*{Appendix: K\"{o}the duality}
\addtocounter{section}{1}
\setcounter{subsection}{0} 
\setcounter{Theorem}{0}

In this section we present a proof of K\"{o}the duality of Musielak-Orlicz spaces. The result is well known. However, to the best of our knowledge, a direct self-contained proof of that result in the general case has never been published. 

We need the following result on Orlicz functions characterizing equality in Young's inequality \cite{MR1909821}. 

\begin{Proposition}
\label{subdiff}
Let $\varphi$ and $\psi$ be a pair of complementary Orlicz functions and $\varphi_{-}^{\prime}$, $\varphi_{+}^{\prime}$ be the left and right derivative of $\varphi$, respectively. Let 
\begin{equation*}
\partial\varphi(u)=\{v\geqslant 0:\varphi(u)+\psi(v)=uv\}\text{, }u\geqslant 0\text{.}
\end{equation*}
Then
\begin{enumerate}[{\rm(i)}]
\item $\partial\varphi(0)=[0,a_{\psi}]=[0,\varphi_{+}^{\prime}(0)]$.
\item If $u\in(0,b_{\varphi})$ then $\partial\varphi(u)=[\varphi_{-}^{\prime}(u),\varphi_{+}^{\prime}(u)]$.
\item If $\varphi_{-}^{\prime}(b_{\varphi})<\infty$ then $\partial\varphi(b_{\varphi})=[\varphi_{-}^{\prime}(b_{\varphi}),\infty)$.
\item If $\varphi(b_{\varphi})=\infty$ then $\partial\varphi(b_{\varphi})=\emptyset$.
\item If $u>b_{\varphi}$ then $\partial\varphi(u)=\emptyset$.
\end{enumerate}
\end{Proposition}

\begin{Lemma}
\label{finitecomp}
Let $\varphi$ and $\psi$ be a pair of complementary Orlicz functions and $\psi^{\prime}_{-}$ be the left derivative of $\psi$. 
If $b_{\varphi}<\infty$ then $\varphi(\psi^{\prime}_{-}(u))<\infty$ for all $u>0$.
\end{Lemma}

\begin{Theorem}
\label{OrliczNormDual}
Let $M$ be a Musielak-Orlicz function. The K\"{o}the dual $(L_M^o)^{\prime}=L_N$.
\begin{proof}
Recall that $(L_M^o)^{\prime}$ is isometrically isomorphic to the space of all order continuous functionals $(L_M^o)^{*}_c$.

Let $g\in L_N$ and $F:L_{M}^o\to\R$ be defined by $F(f)=\int_{\Omega}fg\, d\mu$. By definition of $\|\cdot\|_M^o$, $|F(f)|\leqslant\|f\|_{M}^o\|g\|_N$. Hence $\|F\|\leqslant \|g\|_N$. Thus $F$ is a bounded linear order continuous functional. Next we show the reverse inequality. Without loss of generality we assume that $\|g\|_N=1$. It follows that $\rho_N(g)\leqslant 1$ and $\rho_N((1+\epsilon)g)\geqslant 1$ for all $\epsilon>0$. 
In the sequel, by $N^{\prime}$ we denote the left-side derivative of $N$ with respect to $u$ (we define $N^{\prime}(t,0)=0$, $t\in\Omega$).

We consider two cases. 

Case 1. There is $\epsilon_0>0$ such that $(1+\epsilon_0)|g|\leqslant b_N$ $\mu$-a.e. on $\Omega$. Fix $\epsilon\in(0,\epsilon_0)$. In view of $(1+\epsilon)|g|<b_N$ $\mu$-a.e. on $\Omega$, by Lemma \ref{finitecomp}, the function $M(\cdot,N^{\prime}(\cdot,(1+\epsilon)|g|(\cdot)))$ is nonnegative and finite $\mu$-a.e. on $\Omega$. Since $\mu$ is $\sigma$-finite, there is an ascending sequence of measurable sets with finite and positive measure $ (\Omega_n)_{n=1}^{\infty}$ such that $\Omega=\cup_{n=1}^{\infty}\Omega_n$. Let $T_{n}=\{t\in\Omega_n:M(t,N^{\prime}(t,(1+\epsilon)|g|(t)))\leqslant n\}$, $n\in\N$.  Clearly $ (T_{n})_{n=1}^{\infty}$ is an ascending sequence of measurable sets of finite measure satisfying 
\begin{equation}
\label{supMNg}
\sup_{t\in T_{n}}M(t,N^{\prime}(t,(1+\epsilon)|g|(t)))<\infty\text{, } n\in\N\text{.}
\end{equation} 
Moreover $\mu(\Omega\setminus\cup_nT_{n})=0$. Indeed, for $t\in \Omega\setminus\cup_nT_{n}$ we have that $t\in\Omega_n$ for all $n\geqslant n_0$, for some $n_0\in\N$ and $t\notin \cup_{n=1}^{\infty}T_{n}$. This implies that for all $n\geqslant n_0$, $M(t,N^{\prime}(t,(1+\epsilon)|g|(t)))>n$. Since $M(\cdot,N^{\prime}(\cdot,(1+\epsilon)|g|(\cdot)))$ is finite $\mu$-a.e. on $\Omega$ we conclude that $\mu(\Omega\setminus\cup_nT_{n})=0$. 

Let $\tilde{g}_n$ be a sequence of non-negative simple functions  such that $\tilde{g}_n\uparrow|g|$ $\mu$-a.e. on $\Omega$ and $\mu(\supp \tilde{g}_n)<\infty$, $n\in\N$. Define 
\begin{equation*}
g_n=\tilde{g}_n\chi_{T_{n}}\text{, }n\in\N\text{.}
\end{equation*}
 Clearly $g_n\leqslant|g|$ and $g_n\uparrow|g|$ $\mu$-a.e. on $\Omega$. 
Therefore $\rho_N((1+\epsilon)g)=\lim_{n\to\infty}\rho_N((1+\epsilon)g_n)$ and $\rho_N((1+\epsilon)g_n)\geqslant 1$ for all $n$ large enough. 
 By Proposition \ref{subdiff}, for all $n\in\N$,  
\begin{equation}
\label{Yinq}
\partial N(t,(1+\epsilon)g_n(t))\ne\emptyset\text{ for }\mu\text{-a.a. }t\in\Omega\text{.}
\end{equation}
Moreover, since $g_n\leqslant|g|$, in view of (\ref{supMNg}) we get that for  $n\in\N$,
\begin{equation}
\label{MNgn}
\int_{\Omega}M(t,N^{\prime}(t,(1+\epsilon)g_n(t)))\, d\mu<\infty\text{.}
\end{equation}

 Define
\begin{equation*} 
y_n(t)=N^{\prime}(t,(1+\epsilon)g_n(t))\sign(g(t))\text{, }t\in\Omega\text{, }
\end{equation*} 
and
\begin{equation*}
f_n=\frac{y_n}{\rho_M(y_n)+1}\text{. }
\end{equation*}
The functions $f_n$ are well defined since $\rho_M(y_n)<\infty$ for every $n\in\N$ by (\ref{MNgn}).
By Young's inequality for every $h\in L_0$ and $\mu$-a.a. $t\in\Omega$,
\begin{equation*}
|y_n(t)h(t)|=N^{\prime}(t,(1+\epsilon)g_n(t))|h(t)|\leqslant M(t,N^{\prime}(t,(1+\epsilon)g_n(t)))+N(t,|h(t)|)\text{.}
\end{equation*} 
Hence
\begin{equation*}
\|f_n\|_M^o=\sup\left\{\int_{\Omega}\frac{|y_nh|}{\rho_M(y_n)+1}:\rho_N(h)\leqslant 1\right\}\leqslant 1\text{.}
\end{equation*}

By (\ref{Yinq}) the following equality in Young's inequality holds true for $\mu$-a.a. $t\in\Omega$
\begin{equation*}
\frac{1}{1+\epsilon}N^{\prime}(t,(1+\epsilon)g_n(t))(1+\epsilon)g_n(t)=\frac{1}{1+\epsilon}[M(N^{\prime}(t,(1+\epsilon)g_n(t)))+N(t,(1+\epsilon)g_n(t))]\text{.}
\end{equation*}
It follows that 
\begin{equation*}
\begin{split}
\|F\|&=\sup\left\{\left|\int_{\Omega}fg\, d\mu\right|:\|f\|_M^o\leqslant 1\right\}\geqslant\left|\int_{\Omega}f_ng_n\sign(g)\, d\mu\right|=\int_{\Omega}\frac{N^{\prime}(t,(1+\epsilon)g_n(t))g_n(t)}{\rho_M(y_n)+1}\, d\mu\\
&=\frac{1}{1+\epsilon}\int_{\Omega}\frac{M(t,N^{\prime}(t,(1+\epsilon)g_n(t)))+N(t,(1+\epsilon)g_n(t))}{\rho_M(y)+1}\, d\mu=\frac{1}{1+\epsilon}\frac{\rho_M(y_n)+\rho_N((1+\epsilon)g_n)}{\rho_M(y_n)+1}\text{.}
\end{split}
\end{equation*}
Since $1\leqslant\rho_N((1+\epsilon)g)=\lim_{n\to\infty}\rho_N((1+\epsilon)g_n)$, we conclude that $\|F\|\geqslant 1/(1+\epsilon)=\|g\|_N/(1+\epsilon)$. But we can take $\epsilon$ arbitrarily close to $0$, hence $\|F\|\geqslant\|g\|_N$. 

Case 2. For every $\epsilon>0$ there is a measurable set $E_{\epsilon}$ of positive measure such that 
$(1+\epsilon)|g|>b_N$ $\mu$-a.e. on $E_{\epsilon}$. Let
\begin{equation}
\label{ggr}
A_n=\{t\in\Omega:|g(t)|\geqslant (1-1/n)b_N(t)\}\text{.}
\end{equation}
Clearly $\mu(A_n)>0$ and $b_N<\infty$ $\mu$-a.e. on $A_n$ for every $n\in\N$. There are measurable sets $B_n\subset A_n$ of positive and finite measure such that $\int_{B_n}b_N\, d\mu<\infty$, $n\in\N$. Define
\begin{equation*}
f_n=\left(\int_{B_n}b_N\, d\mu\right)^{-1}\chi_{B_n}\sign(g)\text{, }n\in\N\text{.}
\end{equation*}
Since for any $h\in L_0$ with $\rho_N(h)\leqslant 1$ we have $|h|\leqslant b_N$ $\mu$-a.e. on $\Omega$, so
\begin{equation*}
\|f_n\|^o_M=\sup\left\{\left|\int_{\Omega}f_nh\right|:\rho_N(h)\leqslant 1\right\}\leqslant 1\text{.}
\end{equation*}
From (\ref{ggr}) we get that
\begin{equation*}
\int_{\Omega}f_ng\, d\mu=\left(\int_{B_n}b_N\, d\mu\right)^{-1}\int_{B_n}g\sign(g)\,d\mu=\left(\int_{B_n}b_N\, d\mu\right)^{-1}\int_{B_n}|g|\,d\mu\geqslant 1-1/n\text{.}
\end{equation*}
It follows that $\|F\|\geqslant\|g\|_N$.

The fact that every order continuous functional on $L_M^o$ is of the integral form follows from the Radon-Nikodym Theorem. Indeed, let $F$ be an order continuous functional on $L_M^o$. From the Radon-Nikodym Theorem it follows that there is a measurable function $g$ such that $F(\chi_E)=\int_Eg\, d\mu$ for every measurable set $E$ with $\mu(E)<\infty$ such that $\chi_E\in L_M^o$. Let $f\in L_M^o$ be such that $f\geqslant 0$ $\mu$-a.e. on $\Omega$. There is a sequence $(f_n)$ of simple functions  such that $0\leqslant f_n\leqslant f$ and $f_n\uparrow f$ $\mu$-a.e. on $\Omega$. Since $F$ is order continuous we have that $|F(f-f_n)|\to 0$ as $n\to\infty$. Hence $F(f)=\lim_nF(f_n)=\lim_n \int_{\Omega}f_ng\, d\mu=\int_{\Omega}fg\, d\mu$. Since an arbitrary $f\in L_M^o$ can be written as $f=f^+-f^-$, where $f^+$ and $f^-$ are positive, we see that $F(f)=\int_{\Omega}fg\, d\mu$ for every $f\in L_M^o$. Since $\|f\|_M^o\leqslant 2\|f\|_M$ for every $f\in L_M^o$, 
we have that 
\begin{equation*}
\begin{split}
\|g\|_N^o&=\sup\left\{\left|\int_{\Omega}fg\, d\mu\right|:\|f\|_M\leqslant 1\right\}\leqslant \sup\left\{\left|\int_{\Omega}fg\, d\mu\right|:\left\|\frac{f}{2}\right\|^o_M\leqslant 1\right\}\\
&=2\sup\left\{\left|\int_{\Omega}hg\, d\mu\right|:\|h\|_M^o\leqslant 1\right\}=2\|F\|<\infty\text{.}
\end{split}
\end{equation*}
Hence $g\in L_N$. 
\end{proof}
\end{Theorem}

\begin{Theorem}
\label{dualLM}
Let $M$ be a Musielak-Orlicz function. The K\"{o}the dual $(L_M)^{\prime}=L_N^o$.
\begin{proof}
Let $g\in L_N^o$ and $F:L_{M}\to\R$ be defined by $F(f)=\int_{\Omega}fg\, d\mu$. Clearly $F$ is a bounded linear order continuous functional. Since for $f\in L_M$, $\|f\|_M\leqslant 1$ if and only if $\rho_M(f)\leqslant 1$ we get that $\|F\|=\sup\{|\int_{\Omega}fg\, d\mu|:\|f\|_M\leqslant 1\}=\sup\{|\int_{\Omega}fg\, d\mu|:\rho_M(f)\leqslant 1\}=\|g\|_N^o$. 

The fact that every order continuous functional on $L_M$ is of the integral form follows similarly as in Theorem \ref{OrliczNormDual}.
\end{proof}
\end{Theorem}

\bibliography{DaugavetMO}

\end{document}